\documentclass[11pt,leqno]{article}

\usepackage{amsmath,amssymb,amsthm}
\usepackage{fullpage}
\usepackage{hyperref}
\usepackage{setspace}
\usepackage{dsfont}
\newtheorem{theorem}{Theorem}[section]
\newtheorem{lemma}[theorem]{Lemma}
\newtheorem{proposition}[theorem]{Proposition}
\newtheorem{corollary}[theorem]{Corollary}
\newtheorem{remark}[theorem]{Remark}
\theoremstyle{remark}

\numberwithin{equation}{section}

\newcommand{\ls}{\leqslant}
\newcommand{\gr}{\geqslant}

\newcommand{\ip}[2]{\langle #1,#2\rangle}

\hypersetup{colorlinks=true,linkcolor=blue,citecolor=blue,urlcolor=blue}

\begin{document}
\small

\title{\bf Economic cap coverings for log-concave measures}
\author{Minas Pafis}
\date{}

\maketitle

\begin{abstract}\footnotesize
Let $\mu$ be a log-concave probability measure on $\mathbb R^n$. Write $q_\mu$ for its half-space depth,
$W_t=\{q_\mu\ls t\}$ for its wet part, and $w_\mu(t)=\mu(W_t)$. For sufficiently small $t$, we construct a cover of $W_t$ by half-spaces of 
measure at most $A_nt$, with cardinality bounded by
$$D_n\frac{w_\mu(t)-w_\mu(a_nt)}{t},$$
where the constants depend only on $n$. For each $0<\varepsilon\ls 1$ we also construct a cover by half-spaces of
measure at most $(1+\varepsilon)t$, with cardinality bounded by 
$$B_n\varepsilon^{-(n-1)/2}\frac{w_\mu(t)-w_\mu(a_nt)}{t}.$$
The exponent $(n-1)/2$ is optimal for $n\gr 2$, as shown by the Gaussian example. We further establish containment coverings for half-spaces in every dimension:
there is a family of $O_n(w_\mu(t)/t)$ Borel covering sets, each of $\mu$-measure $O_n(t)$, such that every half-space of measure at most $t$ is contained in one covering set. We also obtain a relative wet-part estimate with sharp exponent $2/(n+1)$ and a reverse doubling inequality.
\end{abstract}

\section{Introduction}\label{sec:introduction}

In this work we extend to log-concave probability measures the economic cap-covering theorem of B\'ar\'any and Larman~\cite{Barany-Larman-1988}
for the uniform measure on a convex body. We begin with this classical result and some of its applications.

We write $A\approx_m B$ if $A$ and $B$ are bounded above and below by positive constant multiples of each other, with constants depending only on $m$.

Let $K$ be a convex body in $\mathbb{R}^n$, and denote by $\mu_K$ the normalized Lebesgue measure on $K$. A cap of $K$ is an intersection 
$C=K\cap H$, where $H$ is a closed half-space. For every $x\in K$ we define
$$v_K(x)=\inf\{\mu_K(K\cap H): H\text{ is a closed half-space containing }x\}.$$
The wet part of $K$ at level $t$ is
$$U_t(K)=\{x\in K:v_K(x) \ls t\},$$
and its complement in $K$ is the interior of the convex floating body of $K$ at level $t$ (see Lemma~\ref{lem:depth-regularity}). 

Macbeath's centrally symmetric regions~\cite{Macbeath-1952}, together with packing arguments, play a central role in this subject. Ewald, Larman, and
Rogers~\cite{Ewald-Larman-Rogers-1970} used related regions to obtain an economic covering theorem for inner parallel bodies. Adapting that method
to floating bodies, B\'ar\'any and Larman~\cite[Theorem~6]{Barany-Larman-1988} proved the following result. 

\begin{theorem}\label{th:barany-larman-economic-cap}
Let $K$ be a convex body in $\mathbb{R}^n$. For every $0<t<t_n$, where $t_n$ depends only
on $n$, there exist caps $C_1,\ldots,C_m$ and disjoint Borel sets $C'_i\subseteq C_i$ such that
\begin{equation}\label{eq:intro-classical-economic}
\bigcup_{i=1}^m C'_i\subseteq U_t(K)\subseteq\bigcup_{i=1}^m C_i,\qquad a_nt\ls\mu_K(C'_i)\ls\mu_K(C_i)\ls A_nt,
\end{equation}
where $a_n,A_n>0$ are constants depending only on $n$. 
\end{theorem}

The fact that $C'_i$ are disjoint implies that
\begin{equation}\label{eq:intro-classical-count}
m\ls \frac{1}{a_nt}\,\mu_K(U_t(K)).
\end{equation}
Conversely, the cardinality of any covering of $U_t(K)$ by caps of relative volume at most $A_nt$ must be at least $\mu_K(U_t(K))/(A_nt)$. 
The method of $M$-regions, its refinements, and other applications are surveyed in~\cite{Barany-2000,Barany-survey-2007}. 

Theorem~\ref{th:barany-larman-economic-cap} has important applications. The original motivation of B\'ar\'any and Larman was the approximation of a
convex body by the convex hull of random points. If $X_1,\ldots,X_N$ are independent and uniformly distributed in $K$ and $P_N=\operatorname{conv}\{X_1,\ldots,X_N\}$, then Theorem~\ref{th:barany-larman-economic-cap} implies that 
\begin{equation}\label{eq:intro-classical-random}
\mathbb{E}\,\mu_K(K\setminus P_N)\approx_n \mu_K(U_{1/N}(K))
\end{equation}
for $N$ sufficiently large in terms of $n$. This result connects a random approximation problem to a deterministic geometric quantity, 
uniformly over all convex bodies. This connection has its origin in the work of Rényi and Sulanke \cite{Renyi-Sulanke-1963} on convex hulls 
of random points in the plane, as well as on Efron's identity \cite{Efron-1965}, which relates the expected number of vertices of a random
polytope to the  expected missing volume. B\'{a}r\'{a}ny \cite{Barany-1989} obtained estimates for the expected number of faces of random 
polytopes. Further results for smooth convex bodies were obtained by B\'{a}r\'{a}ny \cite{Barany-1992}, who obtained asymptotic estimates for the 
expected missing volume and for quermassintegrals. For further background, we refer to the surveys \cite{Barany-2000,Reitzner-2010}.

The geometry of wet parts also shows a dependence on the boundary behavior. For every convex body $K$ in $\mathbb{R}^n$, one has
\begin{equation}\label{eq:intro-classical-wet}
\mu_K(U_t(K))\ls C_n t^{2/(n+1)}
\end{equation}
if $t$ is small enough. The example of a Euclidean ball shows that the exponent $2/(n+1)$ is optimal. On the other hand,
for a fixed $n$-dimensional polytope $P$ one has $\mu_P(U_t(P))\approx_P t(\log(1/t))^{n-1}$ as $t\to 0^+$.
For smooth convex bodies of positive curvature, the asymptotic behavior is the same as for the Euclidean ball and is 
expressed through the affine surface area. For related results on floating bodies, affine surface area, and random polytopes, we refer to \cite{Barany-Larman-1988,Barany-2008-notes,Schutt-Werner-1990,Reitzner-2010,Werner-2022}.

A second major line of applications concerns deterministic approximation by polytopes and algorithms in computational geometry. Classical results
of Dudley \cite{Dudley-1974} and Bronshteyn and Ivanov \cite{Bronshteyn-Ivanov-1975} estimate the number of facets or vertices required
for the approximation of convex bodies by circumscribed polytopes in the Hausdorff distance. 
In a different formulation, Macbeath packings and economic coverings are used to control the complexity of polytope approximations. This viewpoint appears, for example, in the work of Arya, da Fonseca, and Mount on $\varepsilon$-kernels, economic convex coverings, and Banach--Mazur approximation \cite{Arya-daFonseca-Mount-2017,Arya-daFonseca-Mount-2024}, as well as in their area sensitive polytope approximation theorem \cite{Arya-daFonseca-Mount-2026}; see also \cite{Arya-Arya-daFonseca-Mount-2022}. Other applications are related to approximate diameter and width,
nearest neighbor searching, and to the closest vector problem \cite{Arya-daFonseca-Mount-2017,Arya-daFonseca-Mount-2024}.

\medskip

\noindent\textbf{Main results.} We extend the theory of economic cap coverings to log-concave probability measures on $\mathbb{R}^n$. For such a measure $\mu$, the corresponding caps are half-spaces of small $\mu$-measure.

For $x\in\mathbb{R}^n$ define the half-space depth function of $\mu$ by
$$q_\mu(x)=\inf\{\mu(H):H\text{ is a closed half-space containing }x\}.$$
This notion was introduced in Tukey's work on multivariate data analysis \cite{Tukey-1975} and became a central tool in robust statistics; see
Donoho and Gasko~\cite{Donoho-Gasko-1992} and Zuo and Serfling~\cite{Zuo-Serfling-2000}. Its connections with floating
bodies and convex geometry are presented in the survey article of Nagy, Sch\"utt and Werner \cite{Nagy-Schutt-Werner-2019}. 
Weighted floating bodies and approximation with densities on convex bodies have also been investigated by Besau, Ludwig, and
Werner~\cite{Besau-Ludwig-Werner-2018}.  Related depth and random-polytope questions for log-concave measures have been
studied by Brazitikos, Giannopoulos and the author \cite{Brazitikos-Giannopoulos-Pafis-2024a,Brazitikos-Giannopoulos-Pafis-2024b};
see also~\cite{Hayakawa-Lyons-Oberhauser-2023}.

For every $t>0$ we set
$$K_t=\{x:q_\mu(x)\gr t\},\qquad W_t=\{x:q_\mu(x) \ls t\},\qquad w_\mu(t)=\mu(W_t),$$
and write $q_*=\sup_xq_\mu(x)$. If $\mu=\mu_K$ is the uniform probability measure on a convex body $K$, then $q_\mu|_K=v_K$, and $w_\mu(t)=\mu_K(U_t(K))$.
In general, $K_t$ is convex and $W_t$ is the measure theoretic wet part. We call every $0<t<q_*$ an admissible level.

Our starting point is the work of B\'ar\'any, Fradelizi, Goaoc, Hubard and Rote~\cite{Barany-Fradelizi-Goaoc-Hubard-Rote-2020} who 
showed that the classical wet-part lower bound in \eqref{eq:intro-classical-random} continues to hold for arbitrary distributions. 
On the other hand, they gave an example  showing that a logarithmic loss cannot be removed in full generality and asked if the classical two-sided estimate holds true for the class of log-concave probability measures. 

For $0<\varepsilon\ls1$, let $N_\mu(t,\varepsilon)$ be the minimum number of closed half-spaces of $\mu$-measure at most $(1+\varepsilon)t$
needed to cover $W_t$. We shall also use the Macbeath regions
$$M_K(x)=x+((K-x)\cap(x-K)),\qquad M_K^\lambda(x)=x+\lambda((K-x)\cap(x-K))$$
of the convex depth regions $K_t$.

\begin{theorem}\label{th:intro-economic-covering}
There exist $a_n,t_n\in(0,1)$ and $A_n,D_n>0$, depending only on $n$, such that for every $0<t<t_n$ we may choose
half-spaces $C_1,\ldots,C_N$ with
$$W_t\subseteq\bigcup_{i=1}^N \operatorname{int}(C_i)\subseteq\{q_\mu\ls A_nt\},\quad \mu(C_i)\ls A_nt,$$
and
$$N\ls D_n\frac{w_\mu(t)-w_\mu(a_nt)}{t}.$$
\end{theorem}

The proof of Theorem~\ref{th:intro-economic-covering} is given in Section~\ref{sec:halfspace-covering} (see Theorem~\ref{th:economic-cap-covering}
for a detailed version). This shows that the B\'ar\'any--Larman covering principle holds for log-concave measures 
and the covering sets are half-spaces.

\begin{theorem}\label{th:intro-near-minimal}
There exist $a_n,t_n\in(0,1)$ and $B_n>0$, depending only on $n$, such that for $0<t<t_n$ and $0<\varepsilon\ls1$ we may find
half-spaces $C_1,\ldots,C_N$ covering $W_t$ with
$$\mu(C_i)\ls(1+\varepsilon)t,\quad N\ls B_n\varepsilon^{-(n-1)/2}\frac{w_\mu(t)-w_\mu(a_nt)}{t}.$$
For $n\gr2$, the exponent $(n-1)/2$ is optimal uniformly among log-concave probability measures.
\end{theorem}

The proof of Theorem~\ref{th:intro-near-minimal} is given in Section~\ref{sec:halfspace-covering} (see Theorem~\ref{th:near-minimal-measure-covering}).
The case of Gaussian measure shows that the exponent $(n-1)/2$ cannot be improved (see Theorem~\ref{th:gaussian-covering}).

The preceding theorems concern coverings of $W_t$. A different and stronger assertion asks for containment of each small half-space 
in one of a fixed collection of Borel sets of small measure.

\begin{theorem}\label{th:intro-containment}
For every $n\gr1$ there exists a constant $C_n>0$ depending only on $n$ such that, for every $0<t<q_\ast$, one can find Borel sets
$A_1,\ldots,A_J\subseteq\mathbb{R}^n$ with
$$J\ls C_n\frac{w_\mu(t)}{t},\quad \mu(A_i)\ls C_nt,$$
so that every closed half-space $H$ of measure at most $t$ is contained in one of the sets $A_i$. More generally,
if $\mathcal H$ is a family of half-spaces of equal measure $t$, such that $\bigcup_{H\in\mathcal H}H$ is a Borel set, 
we can choose 
$$J\ls C_n\frac{\mu\bigl(\bigcup_{H\in\mathcal H}H\bigr)}{t}.$$
\end{theorem}

The proof of Theorem~\ref{th:intro-containment} is provided in Section~\ref{sec:containment} (see Theorem~\ref{th:containment-family-optimal}). 
We use a representation of a log-concave measure as the image of uniform measure on a compact convex body in $\mathbb{R}^{n+1}$. 
The half-spaces correspond to caps of a convex body, and this allows us to apply the covering lemma of Arya, da Fonseca and Mount~\cite[Lemma~4.3]{Arya-daFonseca-Mount-2026} and then project the resulting covering sets back to $\mathbb{R}^n$.

We also obtain a comparison of wet parts at different depth levels (see Theorem~\ref{th:optimal-wet-relative}).

\begin{theorem}\label{th:intro-wet-relative}
There exist constants $C_n,t_n>0$ depending only on $n$, such that
$$w_\mu(s)\ls C_n\left(\frac{s}{t}\right)^{2/(n+1)}w_\mu(t),\quad 0<s\ls t<t_n.$$
The exponent $2/(n+1)$ is optimal among log-concave probability measures.
\end{theorem}

This recovers in a relative form the characteristic exponent that appears in the case of convex bodies. The sharpness of the
exponent follows from the example of normalized volume on the Euclidean ball.

\medskip

There are two companion papers to this work. In \cite{Giannopoulos-Pafis-doubling} we show that there exists an absolute constant $C>0$ such that,
for any log-concave probability measure $\mu$ on $\mathbb{R}^n$, $n\gr2$,
\begin{equation}\label{eq:intro-companion-doubling}
w_\mu(t)\ls C^n\frac{t}{s}\left(1+\frac{\log(t/s)}{n-1}\right)^{n-1}w_\mu(s), \quad 0<s<t<1.
\end{equation}
In dimension one, the exact formula is $w_\mu(t)=\min\{2t,1\}$ for $t>0$.
Theorem~\ref{th:intro-wet-relative} bounds the ratio $w_\mu(s)/w_\mu(t)$ from above by a power of $s/t$.
So, the two results give complementary two-sided information about the measure of wet parts.
In \cite{Pafis-random} we combine these two types of bounds to study the convex hull
$P_N=\operatorname{conv}\{X_1,\ldots,X_N\}$ of independent random points distributed according to $\mu$, 
and compare the expected number of its faces with $Nw_\mu(1/N)$. 
More precisely, in \cite{Pafis-random} we combine the bound for random hyperplanes in Corollary~\ref{cor:random-hyperplane-bound}, which follows from the covering result in Corollary~\ref{cor:containment-optimal}, with the doubling estimate of \cite{Giannopoulos-Pafis-doubling} to obtain
$$\mathbb{E} f_k(P_N)\approx_{n,k}Nw_\mu(1/N)$$
for $0\ls k\ls n-1$ and all sufficiently large $N$.

\smallskip

The paper is organized as follows. In Section~\ref{sec:depth} we introduce the basic notation and definitions,
and establish the main estimates for half-space depth, that will be used throughout the paper. In Section~\ref{sec:halfspace-covering} 
we prove the half-space covering Theorems~\ref{th:intro-economic-covering} and~\ref{th:intro-near-minimal}, and show that they are 
sharp by considering the Gaussian case. In Section~\ref{sec:containment} we prove Theorem~\ref{th:intro-containment}. Finally, in Section~\ref{sec:shell} 
we prove Theorem~\ref{th:intro-wet-relative} and its consequences for reverse doubling, depth moments, and cardinalities of coverings.

\section{Notation and depth estimates}\label{sec:depth}

We begin with basic notation and definitions. We work in ${\mathbb R}^n$, which is equipped with the standard inner product  $\langle\cdot ,\cdot\rangle $. We denote the corresponding Euclidean norm by $\|\cdot \|_2$, and write $B_2^n$ for the Euclidean unit ball, and $S^{n-1}$ for the unit sphere. Volume in $\mathbb{R}^n$ is denoted by $\operatorname{vol}_n$. Whenever we write $a \approx b$ for positive quantities $a,b$, we mean that there exist absolute constants $c_1, c_2 > 0$ such that $c_1 a \ls b \ls c_2 a$. Numerical constants with a subscript $n$ depend only on the dimension $n$ and may change from line to line.

A convex body in ${\mathbb R}^n$ is a compact convex subset $K$ of ${\mathbb R}^n$ with non-empty interior. We say that $K$ is symmetric if $K=-K$. For a convex body $K$, a point $x \in K$ and every $\lambda>0$ we define  the Macbeath regions
$$M_K(x)=x+((K-x)\cap(x-K)),\qquad M_K^\lambda(x)=x+\lambda((K-x)\cap(x-K)).$$
These are compact convex sets, centrally symmetric about $x$. If $x\in\operatorname{int}(K)$ then they are convex bodies.

If $\mu$ is a Borel probability measure on $\mathbb{R}^n$ and $T$ is an invertible affine transformation, then we denote by $T_{\#}\mu$ the corresponding push-forward measure. A Borel probability measure $\mu$ on $\mathbb R^n$ is called log-concave if it is full-dimensional, i.e.  $\mu(H)<1$ for every hyperplane $H$ in ${\mathbb R}^n$, and
$$\mu(\lambda A+(1-\lambda)B) \gr \mu(A)^{\lambda}\mu(B)^{1-\lambda}$$
for any pair of compact sets $A,B$ in ${\mathbb R}^n$ and any $\lambda \in (0,1)$. Borell \cite{Borell-1974} has proved that, under these assumptions, $\mu $ has a log-concave density $f$. The Brunn-Minkowski inequality implies that if $K$ is a convex body in $\mathbb R^n$ then the uniform probability measure on $K$, denoted by $\mu_K$, is log-concave.

Throughout the article, if $\mu$ is log-concave we will denote by $f$ its density that is log-concave and upper-semicontinuous. On the interior of the convex
support of $\mu$, $f$ is positive and continuous, whereas it vanishes outside the support. Since $f$ has finite positive integral, one can check that there exist constants $a,b>0$ such that $f(x)\ls ae^{-b\|x\|_2}$ for all $x\in {\mathbb R}^n$ (see \cite[Lemma~2.2.1]{BGVV-book}).

If $\mu$ is log-concave and $u\in S^{n-1}$ let
$$S_u(r)=\mu(\{z:\ip{z}{u}\gr r\}).$$ Since every one-dimensional marginal of $\mu$ is log-concave we have that $S_u$ is log-concave. For $0<s<1$, the upper $s$-quantile of $S_u$ is the real number $r_s(u)$ chosen so that $S_u(r_s(u))=s$.

The Tukey half-space depth function of $\mu$ is defined by
$$q_\mu(x)=\inf\{\mu(H): H\text{ is a closed half-space containing }x\},$$
and we write
$$K_t=\{x:q_\mu(x)\gr t\}, \qquad W_t=\{x:q_\mu(x) \ls t\}$$
for the $t$-central region and wet part of $\mu$ at level $t$, respectively. We also set $w_\mu(t)=\mu(W_t)$.

The next results establish some elementary properties of $q_\mu$.

\begin{theorem}\label{th:depth-logconcave}
Let $\mu$ be a log-concave probability measure on $\mathbb{R}^n$. Then $q_\mu$ is log-concave: for all $x,y\in\mathbb{R}^n$ and $\lambda\in (0,1)$,
\begin{equation}\label{eq:depth-logconcavity}
q_\mu((1-\lambda)x+\lambda y)\gr q_\mu(x)^{1-\lambda}q_\mu(y)^\lambda.
\end{equation}
\end{theorem}

\begin{proof}Let $u\in S^{n-1}$ and set
$$F_u(x)=S_u(\ip{x}{u}).$$
The function $F_u$ is log-concave because $S_u$ is log-concave and $x\mapsto\ip{x}{u}$ is affine.  Moreover,
$$q_\mu(x)=\inf_{u\in S^{n-1}}F_u(x).$$
This follows from the observation that for any $u\in S^{n-1}$, among the closed half-spaces of the form $\{z:\ip{z}{u}\gr r\}$ that contain $x$, 
the smallest measure is obtained when $x$ belongs to the boundary of the half-space, and hence allowing all $u$ we consider all half-spaces 
that contain $x$.

Now, for every $u$,
$$F_u((1-\lambda)x+\lambda y)\gr F_u(x)^{1-\lambda}F_u(y)^\lambda\gr q_\mu(x)^{1-\lambda}q_\mu(y)^\lambda.$$
Taking the infimum over $u$ we obtain the assertion.
\end{proof}

\begin{remark}\rm For an arbitrary probability measure, half-space depth is quasi-concave because its superlevel sets are convex. 
\end{remark}

\begin{lemma}\label{lem:depth-regularity}
Let $q_*=\sup_x q_\mu(x)$. The function $q_\mu$ is continuous, the infimum in its definition is attained, and $K_t$ is compact for every $t>0$.  
For every $0<t<q_*$,
\begin{equation}\label{eq:depth-regularity}
 \operatorname{int}(K_t)=\{q_\mu>t\},\quad \mu(\{q_\mu=t\})=0.
\end{equation}
Consequently, for $0<a<1$,
\begin{equation}\label{eq:shell-wet-exact}
 \mu(K_{at}\setminus\operatorname{int}(K_t))=w_\mu(t)-w_\mu(at).
\end{equation}
We call the levels $0<t<q_*$ admissible.
\end{lemma}

\begin{proof}
For $u\in S^{n-1}$ set $F(u,x)=\mu(\{z:\ip{u}{z-x}\gr0\})$.  Since $\mu$ has a density, affine hyperplanes have measure zero, and applying
the dominated convergence theorem we see that $F$ is jointly continuous. Moreover $q_\mu(x)=\min_{u\in S^{n-1}}F(u,x)$, so compactness of the sphere 
shows the continuity of $q_{\mu}$ and that the minimum is attained. If $x\ne0$ and $u=x/\|x\|_2$, then
$$q_\mu(x)\ls \mu(\{z:\langle z,u\rangle\gr\|x\|_2\})\ls \mu(\{z:\|z\|_2\gr\|x\|_2\}),$$
which tends to zero as $\|x\|_2\to\infty$. Thus $q_*$ is attained, and every $K_t$ with $t>0$ is bounded and hence compact.

Fix $0<t<q_*$ and choose $y$ with $q_\mu(y)>t$.  From continuity we have $\{q_\mu>t\}\subseteq\operatorname{int}(K_t)$.  
If $x\in\operatorname{int}(K_t)$ and $q_\mu(x)=t$, then for sufficiently small
$\varepsilon>0$ the point $z=x-\varepsilon(y-x)$ belongs to $K_t$ and
$$x=\frac{1}{1+\varepsilon}z+\frac{\varepsilon}{1+\varepsilon}y.$$
From Theorem~\ref{th:depth-logconcave} we get $q_\mu(x)>t$, which is a contradiction. Thus $\operatorname{int}(K_t)=\{q_\mu>t\}$ and  $\partial K_t=\{q_\mu=t\}$. Since $q_\mu$ is continuous, for every $0<t<q_*$ the convex body $K_t$ has non-empty interior, so its boundary has Lebesgue measure zero. Then, the absolute continuity of $\mu$ proves the
second assertion and equation~\eqref{eq:shell-wet-exact} follows.
\end{proof}

\begin{remark}\rm
If $T$ is an invertible affine transformation and $\nu=T_{\#}\mu$, then $q_\nu(x)=q_\mu(T^{-1}x)$ for every $x \in \mathbb{R}^n$. In particular, $$K_t^\nu=TK_t^\mu, \quad M_{K_t^\nu}^{\lambda}(Tx)=TM_{K_t^\mu}^\lambda(x),$$ for every $\lambda>0,\,t \in (0,q_*)$ and $x \in K_t^\mu$.
\end{remark}

The rest of this section is devoted to proving the next theorem.

\begin{theorem}\label{th:fixed-macbeath-region}For every integer $n\gr1$ there exist constants
$a_n\in (0,1/6),\,c_n\in(0,1)$ and $L_n\gr1$,
depending only on $n$, such that the following holds.  Let $\mu$ be a log-concave probability measure on $\mathbb{R}^n$.
Assume that $q_\mu(x)=t>0$ and $H_x$ is a minimizing half-space at $x$, so $\mu(H_x)=t$.  Then
$$\mu\bigl(H_x\cap M_{K_{a_nt}}^{L_n}(x)\bigr)\gr c_nt.$$
\end{theorem}

We need the following two lemmas.

\begin{lemma}\label{lem:planar-reflection}Fix $\delta\in(0,1)$ and $a\in(0,\delta/3)$.  Let $\rho$ be a log-concave probability measure on an affine subspace of $\mathbb{R}^2$ of dimension
one or two, with coordinates $(s,w)$, and assume that the $s$-marginal is non-degenerate. Assume that
$$q_\rho(0)=\rho\big(\{(s,w):s\gr0\}\big)=t>0$$
and
$$\rho\big(\{(s,w):s\gr0,\ w\ls-1\}\big)\gr\delta t.$$
If
$$0<\varepsilon<\frac{\log(1+(\delta-a)/2)}{\log(1/a)},$$
then
$$\rho\big(\{(s,w):w\gr\varepsilon\}\big)\gr \frac{\delta-a}{2}\,t>at.$$
\end{lemma}

\begin{proof}Let $r$ be the upper $at$-quantile of the $s$-marginal, so that
$$\rho\big(\{(s,w):s\gr r\}\big)=at.$$
Since $\rho\big(\{(s,w):s\gr0\}\big)=t$, we have $r>0$.  Removing the part with $\{(s,w):s\gr r\}$ we get
$$\rho\big(\{(s,w):0\ls s<r,\ w\ls-1\}\big)\gr(\delta-a)t.$$
The half-plane $\{(s,w):s\gr0\}$ is minimizing among half-planes through the origin.  Hence comparison with the rotated half-plane $\{(s,w):s+rw\gr0\}$ gives
$$\rho\big(\{(s,w):s<0,\ s+rw\gr0\}\big)\gr\rho\big(\{(s,w):s\gr0,\ s+rw<0\}\big).$$
The fact that $$\{(s,w):0\ls s<r,\ w\ls-1\} \subseteq \{(s,w):s\gr0,\ s+rw<0\}$$ implies that
$$\rho\big(\{(s,w):s<0,\ w\gr- s/r\}\big)\gr(\delta-a)t.$$
Therefore, for any $\varepsilon>0$ we get that
$$\rho\big(\{(s,w):w\gr\varepsilon\}\big)+\rho\big(\{(s,w):-\varepsilon r<s<0\}\big)\gr(\delta-a)t.$$
Suppose that the first term is smaller than $(\delta-a)t/2$.  If
$$T(x)=\rho\big(\{(s,w):s\gr x\}\big),$$
then
$$T(-\varepsilon r)\gr\left(1+\frac{\delta-a}{2}\right)t.$$
The tail function $T$ is log-concave, while $T(0)=t$ and $T(r)=at$. Then, by log-concavity we get
$$T(-\varepsilon r)\ls t a^{-\varepsilon}.$$
Our choice of $\varepsilon$ makes these two inequalities incompatible.  Hence
$$\rho\big(\{(s,w):w\gr\varepsilon\}\big)\gr\frac{\delta-a}{2}t.$$
Since $a<\delta/3$, the last quantity is strictly larger than $at$.
\end{proof}

\begin{lemma}\label{lem:reflected-depth}
Fix $\delta\in(0,1)$ and $a\in(0,\delta/3)$. Assume that $C=C(a,\delta)$ satisfies
$$C>\frac{\log(1/a)}{\log(1+(\delta-a)/2)}.$$ 
Then the following holds: Let $\mu$ be a log-concave probability measure on $\mathbb{R}^n$ such that
$$q_\mu(0)=\mu(H)=t>0$$
for a half-space $H$ through the origin, and let
$$\nu=\mu(\,\cdot\mid H).$$
Then
$$z\in K_\delta^\nu\quad\Longrightarrow\quad-\frac{z}{C}\in K_{at}^\mu.$$
\end{lemma}

\begin{proof}Choose $\varepsilon=1/C$ in Lemma~\ref{lem:planar-reflection}.  Suppose for a contradiction that $z\in K_\delta^\nu$ but $-z/C\notin K_{at}^\mu$.
Then there exists a linear functional $\ell$ such that:
$$\mu(\{x:\ell(x)\gr\ell(-z/C)\})<at.$$
Since every half-space containing the origin has measure at least $t$, we conclude that $\ell(z)<0$. Set $L=-\ell(z)>0$.  The last
inequality becomes
$$\mu(\{x:\ell(x)\gr L/C\})<at.$$
On the other hand, the half-space $\{x:\ell(x)\ls-L\}$ contains $z$.  Since $q_\nu(z)\gr\delta$,
$$\mu(\{x:x\in H,\ \ell(x)\ls-L\})\gr\delta t.$$
We write $H=\{x:\langle x,v\rangle\gr0\}$ and push $\mu$ forward under
$$x\longmapsto\left(\langle x,v\rangle,\frac{\ell(x)}{L}\right).$$
The resulting measure $\rho$ is log-concave on its affine support, whose dimension is one or two.  Its $s$-marginal is non-degenerate because
$\rho\big(\{(s,w): s\gr0\}\big)=t\in(0,1)$.  Every half-plane through $0$ corresponds to a half-space through $0$, so $q_\rho(0)\gr t$, while
$\rho\big(\{(s,w):s\gr0\}\big)=t$; hence $q_\rho(0)=t$.  Moreover
$$\rho\big(\{(s,w):s\gr0,w\ls-1\}\big)\gr\delta t,\quad\rho\big(\{(s,w):w\gr1/C\}\big)<at.$$
This is a contradiction by Lemma~\ref{lem:planar-reflection}.
\end{proof}

We pass now to the proof of Theorem~\ref{th:fixed-macbeath-region}.

\begin{proof}[Proof of Theorem~$\ref{th:fixed-macbeath-region}$]Translate $x$ to the origin and set $\nu=\mu(\,\cdot\mid H_x)$.  
The conditional measure is again log-concave.  By the lower bound for the
expected half-space depth of a log-concave measure proved in \cite[Theorem~1.3]{Brazitikos-Giannopoulos-Pafis-2024a}, there exists an absolute constant $b>0$ such that
$$\int q_\nu(z)\,d\nu(z)\gr e^{-bn}.$$
Set $\zeta_n=e^{-bn}$ and $\delta_n=\zeta_n/2$.  Since $0\ls q_\nu\ls1$,
$$\nu(K_{\delta_n}^\nu)\gr \zeta_n-\delta_n=\frac{\zeta_n}{2}.$$
Choose $a_n\in(0,\delta_n/3)$ and let $L_n=C(a_n,\delta_n)$ be given by Lemma~\ref{lem:reflected-depth}.  If $z\in K_{\delta_n}^\nu$, then every
half-space containing $z$ has $\mu$-measure at least equal to its measure inside $H_x$, which is at least $\delta_nt$.  Hence
$$z\in K_{\delta_nt}^\mu\subseteq K_{a_nt}^\mu.$$
Since $0\in K_{a_nt}^\mu$ and this set is convex, we also have $z/L_n\in K_{a_nt}^\mu$.  Lemma~\ref{lem:reflected-depth} 
gives $-z/L_n\in K_{a_nt}^\mu$.  Therefore
$$z\in L_n\bigl(K_{a_nt}^\mu\cap(-K_{a_nt}^\mu)\bigr)=M_{K_{a_nt}^\mu}^{L_n}(0).$$
It follows that
$$\mu\bigl(H_x\cap M_{K_{a_nt}}^{L_n}(0)\bigr)=t\,\nu\bigl(M_{K_{a_nt}}^{L_n}(0)\bigr)\gr t\nu(K_{\delta_n}^\nu)\gr \frac{\zeta_n}{2}\,t.$$
Thus we may take $c_n=\zeta_n/2$.
\end{proof}

\begin{remark}\rm Theorem~\ref{th:fixed-macbeath-region} allows a fixed enlargement of the Macbeath region.  In Section~\ref{sec:halfspace-covering} we contract this region using Corollary~\ref{cor:cap-contraction}, which follows from Theorem~\ref{th:cap-density-barrier}.
\end{remark}

\section{Economic half-space coverings}\label{sec:halfspace-covering}

The first main result of this section is the following economic covering theorem for log-concave measures.

\begin{theorem}\label{th:economic-cap-covering}There exist constants $0<a_n,t_n<1$ and $A_n,D_n>0$, depending only on $n$, such that for every log-concave probability measure $\mu$ on $\mathbb{R}^n$ and every $0<t<t_n$, there exist closed half-spaces $C_1,\dots,C_N$ satisfying
\begin{equation}\label{eq:economic-covering}
 W_t\subseteq\bigcup_{i=1}^N {\rm int}(C_i)\subseteq\{q_\mu\ls A_nt\}, \quad \mu(C_i)\ls A_nt,
\end{equation}
and
\begin{equation}\label{eq:economic-covering-number}
 N\ls D_n\frac{\mu\big(K_{a_nt}\setminus\operatorname{int}(K_t)\big)}{t}=D_n\frac{w_\mu(t)-w_\mu(a_nt)}{t}.
\end{equation}
In particular,
$$N\ls D_n\frac{w_\mu(t)}{t}.$$
\end{theorem}

The proof requires several intermediate results.

\begin{lemma}\label{lem:log-concave-max-usc}
Let $f: \mathbb{R}^n \to [0,+\infty)$ be a log-concave and upper-semicontinuous density. If $F\subseteq\mathbb{R}^n$ is non-empty and closed, then $f$ attains its supremum on $F$. Moreover, every marginal of $f$ is also log-concave and upper-semicontinuous.
\end{lemma}

\begin{proof}
Since $f$ is log-concave, there exist $a,b>0$ such that $f(x) \ls a e^{-b\|x\|_2}$ for every $x \in \mathbb{R}^n$. Hence, $f$ is bounded. Set 
$$M_F=\sup_{x \in F} f(x).$$ 
If $M_F=0$, every point of $F$ is a maximizer. We may therefore assume that $M_F>0$. The exponential decay of $f$ implies that $f(x)<M_F/2$ for every $\|x\|_2>R$, for some $R>0$. Therefore,
$$M_F=\sup_{x \in F \cap RB_2^n} f(x)=\max_{x \in F \cap RB_2^n} f(x).$$ 
The last equality follows because an upper-semicontinuous function attains its maximum on a non-empty compact set. Now, for $1 \ls d \ls n-1$, the function 
$$g(x)=\int_{\mathbb{R}^{n-d}}f(x,y)\, dy, \quad x \in \mathbb{R}^d$$ 
is log-concave by Pr\'ekopa's theorem. Upper semicontinuity follows from the reverse Fatou lemma, using the exponential decay of $f$ as an integrable dominating bound.
\end{proof}

\begin{proposition}\label{prop:section-barycenter}Let $\mu$ be a log-concave probability measure on $\mathbb{R}^n$ with density $f$. Let
$x\in\mathbb{R}^n$ and assume that $\mu(H_x)=q_\mu(x)>0$, where $H_x=\{z:\ip{u}{z-x}\gr0\}$, for some $u \in S^{n-1}$. Then
\begin{equation}\label{eq:section-barycenter}
\int_{u^\perp} y f(x+y)\,d\mathcal H^{n-1}(y)=0.
\end{equation}
In other words, $x$ is the $f\,d\mathcal H^{n-1}$-weighted barycenter of the section $x+u^\perp$.
\end{proposition}

\begin{proof} Translate $x$ to the origin and fix $v\in u^\perp$. Let $\varepsilon>0$ and consider the half-space 
$H(\varepsilon)=\{z: \langle z,u+\varepsilon v \rangle \gr 0\}$ that passes through the origin. 
A change of variables shows that the corresponding half-space measure is
$$ M(\varepsilon)=\int_{u^\perp}\int_{-\varepsilon\ip{v}{y}}^\infty f(su+y)\,ds\,dy.$$ 
Since $\varepsilon=0$ minimizes $M$, if the derivative at $0$ exists, then $M'(0)=0$. Set 
$$F(\varepsilon,y)=\frac{1}{\varepsilon} \int_{-\varepsilon \langle v,y \rangle}^{0}f(su+y)\,ds,\quad \varepsilon \neq 0, \ y \in u^\perp.$$ 
Let $K =\operatorname{supp}(\mu)$. Since $q_\mu(0)>0$, we have that $0\in\operatorname{int}(K)$, which implies that
$\partial_{u^\perp}(K\cap u^\perp)
=\partial K\cap u^\perp$.
Hence $K\cap u^\perp$ is a convex set with non-empty relative
interior and thus its relative boundary has $\mathcal H^{n-1}$-measure
zero. Therefore, as $f$ is continuous outside $\partial K$,
its restriction to $u^\perp$ is continuous
at almost every point of the hyperplane. As a result, 
$$F(\varepsilon,y) \to \langle v,y\rangle f(y), \quad \varepsilon \to 0 ,$$ 
for almost every $y \in u^\perp$. Since $f$ is a log-concave density, there exist constants $A,b>0$ such that $f(z)\ls Ae^{-b\|z\|_2}$.  If $\varepsilon$ is non-zero,
then $|F(\varepsilon,y)|$ is dominated by $C\|y\|_2e^{-b\|y\|_2}$, which is integrable on $u^\perp$. Applying the dominated convergence theorem we get
$$ M'(0)=\int_{u^\perp}\ip{v}{y}f(y)\,dy.$$
Since this derivative vanishes for every $v\in u^\perp$, we obtain \eqref{eq:section-barycenter}.
\end{proof}

\begin{lemma}\label{lem:tail-peak}Let $h$ be a log-concave probability density on $\mathbb{R}$ and assume that
$$t=\int_0^\infty h(s)\,ds\in(0,1).$$
Then
\begin{equation}\label{eq:tail-peak}
 \max_{s\gr0}h(s)\ls \frac{h(0)}{1-t}.
\end{equation}
Thus, for $0\ls\lambda\ls1$ and $s\gr0$,
\begin{equation}\label{eq:normal-contraction}
 h(\lambda s)\gr (1-t)^{1-\lambda}h(s).
\end{equation}
\end{lemma}

\begin{proof}We set $h_0=h(0)$. We have that $h$ is a log-concave density, and hence $h$ is unimodal. If $h$ is non-increasing on 
$[0,\infty)$ or has its maximum at $0$ there is nothing to prove. Otherwise, there exists $s_0>0$ with $$\max_{s \gr 0}h(s)=h(s_0)>h_0.$$ Set
$R=h(s_0)/h_0>1$ and $a=(\log R)/s_0$. By log-concavity we get
$$h(s)\ls h_0e^{as}\quad(s\ls0),\quad h(s)\gr h_0e^{as}\quad(0\ls s\ls s_0).$$
Hence
$$1-t\ls \frac{h_0}{a},\quad t\gr \frac{h_0}{a}(R-1),$$
so $R\ls(1-t)^{-1}$, which is \eqref{eq:tail-peak}. Finally, by log-concavity we have that
$$h(\lambda s)\gr h(0)^{1-\lambda}h(s)^\lambda\gr(1-t)^{1-\lambda}h(s).$$
This completes the proof.\end{proof}

\begin{theorem}\label{th:cap-density-barrier}
Let $\mu$ be a log-concave probability measure on $\mathbb{R}^n$ with density $f$. Assume that $q_\mu(x)=t>0$ and that
$$H_x=\{z:\ip{u}{z-x}\gr0\},\quad \mu(H_x)=t,$$
is minimizing at $x$. Then $t\ls1/2$ and
\begin{equation}\label{eq:cap-density-barrier}
\max_{z\in H_x}f(z)\ls R_nf(x), \quad R_1=2,\quad R_n=4^n\ \ (n\gr2).
\end{equation}
\end{theorem}

\begin{proof}Translate $x$ to the origin and rotate so that $u=e_1$. Since hyperplanes have $\mu$-measure zero, 
opposite half-spaces through the origin have measures summing to $1$, and therefore $t\ls1/2$.

For $n=1$, Lemma~\ref{lem:tail-peak} gives the claim. Assume that $n\gr2$ and set $d=n-1$. For $s\in\mathbb{R}$, $y\in\mathbb{R}^d$, we write
$$g_s(y)=f(s,y),\quad h(s)=\int_{\mathbb{R}^d}g_s(y)\,dy,\quad m(s)=\max_{y\in\mathbb{R}^d}g_s(y).$$
By Pr\'ekopa's theorem, $h$ is a log-concave density and $\int_0^\infty h=t$. From Proposition~\ref{prop:section-barycenter} we know 
that the barycenter of $g_0$ is zero. Fix $s\gr0$ and $0<\theta<1$. Choosing $y_s$ to be a maximizer of $g_s$ and using log-concavity, we get
$$ g_{\theta s}((1-\theta)z+\theta y_s) \gr g_0(z)^{1-\theta}g_s(y_s)^\theta.$$
After integration and an affine change of variables,
$$h(\theta s)\gr(1-\theta)^d g_s(y_s)^\theta\int_{\mathbb{R}^d}g_0(z)^{1-\theta}\,dz.$$
The function $p=g_0/h(0)$ is a log-concave probability density on $\mathbb R^d$ with barycenter $0$. Therefore, $p(y)\ls Ae^{-b\|y\|_2}$ for some
$A,b>0$. In particular $p$ is bounded, has finite first moment, and
$\int p^r<\infty$ for every $r>0$. Applying Jensen's inequality  to the convex function $p^{-\theta}$, defined on the convex set $\{p>0\}$, we get
\begin{equation}\label{eq:section-renyi}
 \int_{\mathbb R^d}p^{1-\theta}\,dy
 \gr p(0)^{-\theta}.
\end{equation}
Equivalently,
$$\int_{\mathbb R^d}g_0^{1-\theta}\gr h(0)f(0)^{-\theta}.$$
Therefore, we see that
$$\left(\frac{m(s)}{f(0)}\right)^\theta\ls \frac{h(\theta s)}{h(0)}(1-\theta)^{-d}.$$
Lemma~\ref{lem:tail-peak} bounds the first ratio on the right by
$(1-t)^{-1}\ls2$. Taking $\theta=1/2$ yields
$$\left(\frac{m(s)}{f(0)}\right)^{1/2}\ls 2^{d+1}=2^n.$$
Thus $m(s)\ls 4^nf(0)$ for every $s\gr 0$, as required.
\end{proof}

\begin{corollary}\label{cor:cap-contraction} Let $\mu$ be a log-concave probability measure on $\mathbb{R}^n$. If $q_\mu(x)=t>0$ and $H_x$ is a minimizing half-space, then for every Borel $E\subseteq H_x$ and every $0<\lambda\ls1$, 
\begin{equation}\label{eq:cap-contraction}
\mu(x+\lambda(E-x))\gr \lambda^nR_n^{-(1-\lambda)}\mu(E).\end{equation}
\end{corollary}

\begin{proof}For $z\in H_x$, by log-concavity and Theorem~\ref{th:cap-density-barrier} we get
$$f(x+\lambda(z-x))\gr f(x)^{1-\lambda}f(z)^\lambda\gr R_n^{-(1-\lambda)}f(z).$$
Then we integrate over $E$ and use the Jacobian $\lambda^n$.
\end{proof}

\begin{theorem}\label{th:contracted-macbeath-region}Let $\mu$ be a log-concave probability measure on $\mathbb{R}^n$. Let $a_n,c_n,L_n$ be the constants from Theorem~$\ref{th:fixed-macbeath-region}$. For
every fixed $0<\lambda\ls1$ there exists a constant $c_{n,\lambda}>0$ such that, whenever $q_\mu(x)=t>0$ and $H_x$ is a 
minimizing $t$-half-space through $x$,
\begin{equation}\label{eq:contracted-macbeath-region}
\mu\bigl(H_x\cap M_{K_{a_nt}}^\lambda(x)\bigr) \gr c_{n,\lambda}t.
\end{equation}
\end{theorem}

\begin{proof}Set $K=K_{a_nt}$ and $B_x=H_x\cap M_K^{L_n}(x)$. From Theorem~\ref{th:fixed-macbeath-region} we know that 
$\mu(B_x)\gr c_nt$. Set $\eta=\lambda/L_n$. Since $H_x$ is convex with $x$ on its boundary, and by the definition of Macbeath regions, we can see that
$$x+\eta(B_x-x)\subseteq H_x \cap M_K^\lambda(x).$$ Corollary~\ref{cor:cap-contraction} shows that
$$\mu(H_x\cap M_K^\lambda(x))\gr \eta^nR_n^{-(1-\eta)}\mu(B_x)\gr c_nt\left(\frac{\lambda}{L_n}\right)^n R_n^{-(1-\lambda/L_n)},$$
and the theorem follows.\end{proof}

\begin{lemma}\label{lem:min-cap-supports-depth-body}Let $\mu$ be a log-concave probability measure on $\mathbb{R}^n$. Assume that $q_\mu(x)=t>0$ and
$H_x=\{z:\ip{u}{z}\gr r\}$, where $r=\ip{u}{x}$ is a half-space with $\mu(H_x)=t$. Then
\begin{equation}\label{eq:min-cap-supports-depth-body}
 K_t\subseteq\{z:\ip{u}{z}\ls r\}, \quad H_x\cap\operatorname{int}(K_t)=\varnothing.
\end{equation}
\end{lemma}

\begin{proof}Let $T(s)=\mu(\{z:\ip{u}{z}\gr s\})$. Since $q_\mu(x)>0$, we see that $r$ lies in the interior of the support of
$\ip{u}{X}$. Note that $T$ is strictly decreasing in this interval. This implies that if $y$ belongs to the support of $\mu$
and $\ip{u}{y}>r$ then
$$q_\mu(y)\ls T(\ip{u}{y})<T(r)=t.$$
Outside the support, $q_\mu(y)=0$. It follows that $K_t\subseteq\{z:\ip{u}{z}\ls r\}$, and $K_t$ intersects with $H_x$ only
on the hyperplane $\{z:\ip{u}{z}= r\}$. This proves the lemma.
\end{proof}

The following lemma is from \cite[Lemma~2.4]{Arya-daFonseca-Mount-2017}. We provide a proof for completeness.

\begin{lemma}\label{lem:macbeath-expansion}Let $K$ be a convex body, $x,y\in K$, and $0<\lambda\ls1/5$. If
$$M_K^\lambda(x)\cap M_K^\lambda(y)\neq\varnothing,$$
then
\begin{equation}\label{eq:macbeath-expansion}
M_K^\lambda(y)\subseteq M_K^{4\lambda}(x).
\end{equation}
In particular $y\in M_K^{4\lambda}(x)$.
\end{lemma}

\begin{proof}We define
$$A_p:=(K-p)\cap(p-K),\quad M_K^\eta(p)=p+\eta A_p.$$
We first prove that if $q\in M_K^\eta(p) $ for some $0\ls\eta<1$, then
\begin{equation}\label{eq:macbeath-center-comparison}
 (1-\eta)A_p\subseteq A_q\subseteq(1+\eta)A_p.
\end{equation}
Indeed, assume that $q=p+\eta a$ for some $a \in A_p$. Then, for every $v\in A_p$, the points $q\pm(1-\eta)v$ are convex combinations of $p+a$ and $p\pm v$, 
and this proves the first inclusion. Conversely, if $w\in A_q$, then
$$p\pm\frac{w}{1+\eta}=\frac{1}{1+\eta}(q\pm w)+\frac{\eta}{1+\eta}(p-a)\in K,$$
so $w\in(1+\eta)A_p$.

Now choose $z\in M_K^\lambda(x)\cap M_K^\lambda(y)$. Since $z-x\in\lambda A_x$ and $A_y$ is symmetric, 
applying \eqref{eq:macbeath-center-comparison} first to $z\in M_K^\lambda(y)$, and then to $z \in M_K^\lambda(x)$, we obtain
$$y-x=(z-x)+(y-z)\in\lambda A_x+\frac{\lambda}{1-\lambda}A_z\subseteq\frac{2\lambda}{1-\lambda}A_x.$$
Define $\alpha=2\lambda/(1-\lambda)$. Since $\lambda\ls1/5$, we have $\alpha\ls 1/2$ and hence $y\in M_K^\alpha(x)$. Applying
\eqref{eq:macbeath-center-comparison} once more we get $A_y\subseteq(1+\alpha)A_x$. Therefore
$$M_K^\lambda(y)\subseteq x+\bigl(\alpha+\lambda(1+\alpha)\bigr)A_x=x+\lambda\frac{3+\lambda}{1-\lambda}A_x\subseteq x+4\lambda A_x=M_K^{4\lambda}(x),$$
because $5\lambda\ls1$. This proves the lemma.
\end{proof}

\begin{proposition}\label{prop:cap-expansion-measure}Let $\mu$ be a log-concave probability measure on $\mathbb{R}^n$. Fix $0<a<1$. Let $q_\mu(x)=t>0$ and
$$H_x=\{z:\ip{u}{z}\gr r\},\quad \mu(H_x)=t,$$
be minimizing at $x$. Set $K=K_{at}$ and $h=h_K(u)$. For $\beta\gr1$ define
$$H_x^{[\beta]}:=\{z:\ip{u}{z}\gr h-\beta(h-r)\}.$$
Then
\begin{equation}\label{eq:cap-expansion-measure}
 \mu(H_x^{[\beta]})\ls a^{-(\beta-1)}t.
\end{equation}
\end{proposition}

\begin{proof}Let $T(s)=\mu(\{z:\ip{u}{z}\gr s\})$ and let $r_a$ be the upper $at$-quantile, so $T(r_a)=at$. Since $K_{at}\subseteq\{z:\ip{u}{z}\ls r_a\}$, we have
$h\ls r_a$. Writing $\varepsilon=\beta-1$ and $s_\beta=h-\beta(h-r)=r-\varepsilon(h-r)$ we get
$$s_\beta\gr r-\varepsilon(r_a-r).$$
The tail function $T$ of a one-dimensional log-concave density is log-concave. Therefore, extrapolation of the secant from $r$
to $r_a$ gives
$$ T(r-\varepsilon(r_a-r))\ls a^{-\varepsilon}t.$$
Since $T$ is decreasing, \eqref{eq:cap-expansion-measure} follows.
\end{proof}

\begin{lemma}\label{lem:boundary-to-wet-part}
Let $K\subseteq\mathbb{R}^n$ be a convex body and $U_1,\ldots,U_N$ open
half-spaces such that $\partial K\subseteq\bigcup_i U_i$. If
$$P=\bigcap_{i=1}^N(\mathbb{R}^n\setminus U_i)$$
intersects $\operatorname{int}(K)$, then $P\subseteq\operatorname{int}(K)$.
Therefore, $\mathbb{R}^n\setminus K\subseteq\bigcup_i U_i$.
\end{lemma}

\begin{proof}
The set $P$ is convex and $P\cap (\partial K)=\varnothing$. Choose
$b\in P\cap\operatorname{int}(K)$. If $z\in P\setminus \operatorname{int}(K)$ then the segment $[b,z]\subseteq P$
intersects $\partial K$, and we arrive at a contradiction.
\end{proof}

We are now able to prove the economic  covering theorem.

\begin{proof}[Proof of Theorem~$\ref{th:economic-cap-covering}$]
We use the constant $a_n$ from Theorem~\ref{th:contracted-macbeath-region} and fix $\lambda=1/5$. For every $x\in\partial K_t$, 
we choose a half-space $H_x$ such that $x\in H_x$ and $\mu(H_x)=t$, and set
$$D_x:=H_x\cap M_{K_{a_nt}}^{1/5}(x).$$
From Theorem~\ref{th:contracted-macbeath-region} we get $\mu(D_x)\gr c'_nt$ and Lemma~\ref{lem:min-cap-supports-depth-body} shows that
$$D_x\subseteq K_{a_nt}\setminus\operatorname{int}(K_t).$$
Choose a maximal disjoint family $D_{x_1},\dots,D_{x_N}$. This family is finite because every set $D_x$ has measure at least $c'_nt$.  
Then
$$Nc'_nt\ls\mu\big(K_{a_nt}\setminus\operatorname{int}(K_t)\big),$$
which gives the cardinality estimate.

For $y\in\partial K_t$, maximality gives $D_y\cap D_{x_i}\ne\varnothing$ for some $i$. Hence the corresponding
$1/5$-Macbeath regions intersect, and Lemma~\ref{lem:macbeath-expansion} gives
$$y\in M_{K_{a_nt}}^{4/5}(x_i).$$
Write
$$H_{x_i}=\{z:\ip{u_i}{z}\gr r_i\},\quad h_i=h_{K_{a_nt}}(u_i).$$
As $y\in M_{K_{a_nt}}^{4/5}(x_i)$, we get that
$$\ip{u_i}{y}\gr r_i-\frac{4}{5}(h_i-r_i).$$
Since $q_\mu(x_i)=t>a_nt$, from Lemma~\ref{lem:depth-regularity} we see that $x_i\in\operatorname{int}(K_{a_nt})$ and hence $h_i-r_i>0$.
Set $\widetilde\gamma=9/10$ and define
$$C_i:=H_{x_i}^{[1+\widetilde\gamma]}.$$
Then $\partial K_t\subseteq\bigcup_i\operatorname{int}(C_i)$, while Proposition~\ref{prop:cap-expansion-measure} gives
$$\mu(C_i)\ls a_n^{-9/10}t.$$
Finally set
$$P:=\bigcap_i\big(\mathbb{R}^n\setminus\operatorname{int}(C_i)\big), \quad A_n=a_n^{-9/10}.$$
Choose a constant $B_n>A_n$. If $b_\mu$ is the barycenter of $\mu$, Gr\"unbaum's inequality gives
$q_\mu(b_\mu)\gr e^{-1}$ (see, for example, \cite{Fradelizi-Langharst-Liu-MarinSola-Tang-2026}). Hence, after fixing $t_n<1/(eB_n)$, the body $K_{B_nt}$ is non-empty for every
$0<t<t_n$. It is contained in $P$, since every point of $\operatorname{int}(C_i)$ has depth $<B_nt$. Thus $P$ contains a point
of $\operatorname{int}(K_t)$, while $P\cap\partial K_t=\varnothing$. Lemma~\ref{lem:boundary-to-wet-part}, applied with $K=K_t$,
therefore gives $P\subseteq\operatorname{int}(K_t)$ and
$W_t\subseteq\bigcup_i {\rm int}(C_i)$.  Since each $C_i$ is a half-space of measure at
most $A_nt$, every $z\in C_i$ satisfies $q_\mu(z)\ls A_nt$.  This gives the second inclusion in \eqref{eq:economic-covering}.

The equality in \eqref{eq:economic-covering-number} follows from Lemma~\ref{lem:depth-regularity}.
\end{proof}

\begin{remark}\label{rem:Dn-constant-estimate}\rm
With the notation of Theorem~\ref{th:fixed-macbeath-region} set $$\zeta_n=e^{-bn},\quad\delta_n=\zeta_n/2,\quad a_n=\zeta_n/12,\quad
L_n=\frac{2\log(1/a_n)}{\log(1+(\delta_n-a_n)/2)}, \quad c_n=\zeta_n/2.$$ 
The lower bound in Theorem~\ref{th:contracted-macbeath-region} with $\lambda=1/5$ yields the constant $$D_n=\frac{R_n}{c_n}(5L_n)^n.$$
Since $R_n\ls4^n$ and $L_n \ls e^{dn}$ for some $d>0$, we get that $$D_n\ls e^{cn^2}$$ for some absolute constant $c>0$. 
\end{remark}

\begin{corollary}\label{cor:covering-reverse-doubling}
Let $\mu$ be a log-concave probability measure on $\mathbb{R}^n$. With the constants of Theorem~$\ref{th:economic-cap-covering}$, there exists
$c_n>0$ such that, for $0<t<t_n$,
$$w_\mu(a_nt)\ls(1-c_n)w_\mu(t).$$
\end{corollary}

\begin{proof}The covering and its measure and cardinality bounds give
$$w_\mu(t)\ls\sum_{i=1}^N\mu(C_i)\ls A_ntN\ls A_nD_n\bigl(w_\mu(t)-w_\mu(a_nt)\bigr).$$
Rearranging we conclude the proof, with $c_n=\min\{1/2,(A_nD_n)^{-1}\}$.
\end{proof}

\begin{lemma}\label{lem:local-normal-net}
Let $B$ be a convex body in $\mathbb{R}^n$, let $R\gr1$, and let $S\subseteq\partial B\cap[-R,R]^n$. For every $y\in S$ choose $b_y$ to be an inward unit
supporting normal to $B$. Then, for every $0<\delta\ls1$, there exist points $y^1,\ldots,y^m\in S$ such that
$$m\ls Q_n(R)\delta^{-(n-1)},\quad Q_n(R):=2n\bigl(1+2n(R+1)\bigr)^{n-1},$$
and for every $y\in S$ some $k$ satisfies
$$0\ls\langle b_{y^k},y-y^k\rangle\ls\delta^2/2.$$
\end{lemma}

\begin{proof}Set $v_y=-b_y$ and $p_y=y+v_y$. Each $p_y$ belongs to
$\partial(B+B_2^n)$ with outward supporting normal $v_y$, and
$p_y\in[-R-1,R+1]^n$. We partition the points $p_y$ into at most $2n$ groups,
choosing a coordinate $i$ and a sign for which $|v_{y,i}|\gr n^{-1/2}$.
Within each group the supporting inequalities at $p_y,p_z$ imply that 
\begin{equation}\label{eq:first-partition-inequality}
|p_{y,i}-p_{z,i}|\ls\sqrt{n-1}\|(p_y-p_z)_{\widehat i}\|_2.
\end{equation}
To see this, note that the two supporting inequalities bound the signed $i$-th coordinate
difference from above and below, and
$\|(v_y)_{\widehat i}\|_2/|v_{y,i}|\ls\sqrt{n-1}$, and likewise for $z$.
We partition the remaining coordinates into cells of side $\delta/n$ and
choose one representative from every occupied cell in every group.
There are at most $Q_n(R)\delta^{-(n-1)}$ representatives. From \eqref{eq:first-partition-inequality}
we see that any two points in the same cell have distance $\|p_y-p_z\|_2\ls\delta$.

For a representative $z$ associated with $y$, the supporting inequalities
for $B$ give $\langle v_z,y-z\rangle\ls0$ and $\langle v_y,y-z\rangle\gr0$. Hence
\begin{align*}
0\ls \langle b_z,y-z\rangle
&\ls\langle v_y-v_z,y-z\rangle\\
&\ls\tfrac12\bigl(\|y-z\|_2^2+\|v_y-v_z\|_2^2\bigr)
\ls\tfrac12\|p_y-p_z\|_2^2\ls\delta^2/2.
\end{align*}
For the penultimate inequality we used
$\langle v_y-v_z,y-z\rangle\gr0$. This proves the lemma.
\end{proof}

The next normalization lemma is a direct consequence of John's ellipsoid theorem (see \cite[Theorem~1.5.7]{BGVV-book}).

\begin{lemma}\label{lem:symmetric-normalization}
If $A$ is a symmetric convex body in $\mathbb{R}^n$, there exists an invertible
linear map $T$ such that
$$B_2^n\subseteq T A\subseteq \sqrt{n}B_2^n.$$
\end{lemma}

We can now prove the near-minimal covering theorem.

\begin{theorem}\label{th:near-minimal-measure-covering}
Let $a_n,c_n,L_n$ be the constants of Theorem~$\ref{th:fixed-macbeath-region}$.
There exists a constant $B_n>0$, depending only on $n$, such that for
$0<t<1/(4e)$ and $0<\varepsilon\ls1$ there are finitely many closed
half-spaces $C_1,\ldots,C_N$ with
\begin{equation}\label{eq:near-minimal-mass-covering}
W_t\subseteq\bigcup_{j=1}^N {\rm int}(C_j)\subseteq\{q_\mu\ls(1+\varepsilon)t\},
\qquad \mu(C_j)\ls(1+\varepsilon)t,
\end{equation}
and
\begin{equation}\label{eq:near-minimal-mass-count}
N\ls B_n\varepsilon^{-(n-1)/2}
\,\frac{w_\mu(t)-w_\mu(a_nt)}{t}.
\end{equation}
In fact the measures of the covering half-spaces are at most $\sqrt{1+\varepsilon}\,t$.
For $n\gr2$ the exponent $(n-1)/2$ cannot be decreased uniformly over log-concave measures.
\end{theorem}

\begin{proof}
Write $a=a_n$, $K=K_{at}$ and $A_x=(K-x)\cap(x-K)$.
Every $x\in\partial K_t$ lies in $\operatorname{int}(K)$, since $t>at$,
so $A_x$ is a symmetric convex body.
For $x\in\partial K_t$, choose a minimizing half-space $H_x$ and set
$D_x=H_x\cap M_K^{1/5}(x)$. Theorem~\ref{th:contracted-macbeath-region} gives
$$\mu(D_x)\gr c_n^\prime t.$$
Each $D_x$ lies in $K_{at}\setminus\operatorname{int}(K_t)$.
As in the proof of Theorem~\ref{th:economic-cap-covering}, we choose a maximal pairwise disjoint
family $D_{x_1},\ldots,D_{x_J}$. This family has to be finite, and
\begin{equation}\label{eq:fixed-patch-count}
J\ls D_n\frac{w_\mu(t)-w_\mu(at)}{t}.
\end{equation}
Since the family is maximal, from Lemma~\ref{lem:macbeath-expansion} we get
$$\partial K_t\subseteq\bigcup_{i=1}^J\bigl(x_i+\tfrac45 A_{x_i}\bigr).$$
Fix $i$. Using Lemma~\ref{lem:symmetric-normalization} we choose an affine map
$F_i(z)=T_i(z-x_i)$ so that
$$B_2^n\subseteq A^{(i)}:=T_i A_{x_i}\subseteq\sqrt n\,B_2^n.$$ Write $\mu_i=(F_i)_\#\mu$, $K^{(i)}=F_iK$ and $B^{(i)}=F_iK_t$. Then $A^{(i)} \subseteq K^{(i)}$ and the set
$$S_i=F_i\bigl(\partial K_t\cap(x_i+\tfrac45 A_{x_i})\bigr)$$
is contained in $\partial B^{(i)}\cap[-\sqrt{n},\sqrt{n}]^n$.
For each $y\in S_i$, choose a minimizing half-space, with respect to $\mu_i$, and write it in the form 
$$H_y=\{z:\langle u_y,z\rangle\gr r_y\},\quad
\|u_y\|_2=1,\quad r_y=\langle u_y,y\rangle.$$
Since half-space depth is preserved by invertible affine maps, we have 
$K^{(i)}=K_{at}^{\mu_i}$ and $B^{(i)}=K_t^{\mu_i}$.
Its inward normal $b_y=-u_y$ supports $B^{(i)}$. Since
$y\in(4/5)A^{(i)}$ and $B_2^n \subseteq A^{(i)}$, its height relative to $K^{(i)}$ satisfies
\begin{equation}\label{eq:normalized-cap-margin}
d_y:=h_{K^{(i)}}(u_y)-r_y
\gr h_{A^{(i)}}(u_y)-\tfrac45h_{A^{(i)}}(u_y)
\gr\tfrac15.
\end{equation}
Choose
$$\delta^2=\frac{\log(1+\varepsilon)}{10\log(1/a)}.$$
Since $a<1/6$ and $\varepsilon\ls1$, we have $0<\delta<1$.
From Lemma~\ref{lem:local-normal-net}, with $R=\sqrt{n}$, we select at most
$Q_n(\sqrt{n})\delta^{-(n-1)}$ points of $S_i$. At every selected point $z$
we replace $H_z$ by the enlarged
$$\widetilde C_z=\{v:\langle u_z,v\rangle\gr r_z-\delta^2\}.$$
For each $y\in S_i$, a selected $z$ satisfies
$0\ls r_z-\langle u_z,y\rangle\ls\delta^2/2$, so $y$ lies in the
interior of $\widetilde C_z$.
For this family we can apply Proposition~\ref{prop:cap-expansion-measure} with $\beta=1+\delta^2/d_z$. Therefore,
\eqref{eq:normalized-cap-margin} gives
$$\mu_i(\widetilde C_z)\ls a^{-\delta^2/d_z}t\ls a^{-5\delta^2}t=\sqrt{1+\varepsilon}\,t.$$
We consider the inverse images of these half-spaces by $F_i$ and repeat the same process over all normalized boundary regions.
Their interiors cover $\partial K_t$, and their number is at most
$$JQ_n(\sqrt{n})\delta^{-(n-1)}\ls JQ_n(\sqrt{n})\bigl(20\log(1/a)\bigr)^{(n-1)/2}\varepsilon^{-(n-1)/2},$$
if we take into account $\log(1+\varepsilon)\gr\varepsilon/2$. Together with
\eqref{eq:fixed-patch-count}, this proves the stated cardinality bound.

Finally, set $P=\bigcap_j\big(\mathbb{R}^n\setminus\operatorname{int}(C_j)\big)$.
Every ${\rm int}(C_j)$  has measure less than $2t$, so $K_{4t}\subseteq P$.
Since $t<1/(4e)$, we can find a point in $K_{4t}$,
which lies in $\operatorname{int}(K_t)$. On the other hand
$P\cap\partial K_t=\varnothing$. From Lemma~\ref{lem:boundary-to-wet-part} we see that 
$P\subseteq\operatorname{int}(K_t)$, and this implies that the open half-spaces cover $W_t$.
The right-hand side inclusion of the theorem follows from the fact that $\mu(C_j)\ls (1+\varepsilon)t$.
The sharpness assertion follows from the Gaussian case (see Theorem~\ref{th:gaussian-covering}).
\end{proof}

\begin{remark}\label{rem:near-minimal-constant}\rm
If $D_n$ is the constant from Theorem~\ref{th:economic-cap-covering} then $$B_n= D_nQ_n(\sqrt{n}) \bigl(20\log(1/a_n)\bigr)^{(n-1)/2}.$$ Considering Remark~\ref{rem:Dn-constant-estimate} along with the fact that $Q_n(\sqrt{n}) \ls n^{dn}$ for some $d>0$, we get that $$B_n \ls e^{Cn^2},$$ where $C>0$ is an absolute constant.
\end{remark}

The near-minimal bound in Theorem~\ref{th:near-minimal-measure-covering}, together with an elementary covering lower bound, identify the minimum number of near-minimal half-spaces $N_\mu(t,\varepsilon)$, up to constants depending on the dimension 
and the relative measure parameter $\varepsilon$, with the wet-part measure divided by the depth
scale.

\begin{corollary}\label{cor:entropy-equivalence}
There exists a constant $B_n>0$ depending only on $n$, such that, for every log-concave probability measure $\mu$, 
and any $0<t<1/(4e)$, and $0<\varepsilon\ls1$,
\begin{equation}\label{eq:entropy-equivalence}
\frac{w_\mu(t)}{(1+\varepsilon)t}\ls N_\mu(t,\varepsilon)\ls B_n\varepsilon^{-(n-1)/2}\frac{w_\mu(t)}{t}.
\end{equation}
In particular, for every $\varepsilon \in (0,1]$ we have $N_\mu(t,\varepsilon)\approx_{n,\varepsilon}w_\mu(t)/t$.
\end{corollary}

\begin{proof}Consider a cover $\{H_1,\ldots ,H_J\}$ by half-spaces, each of measure at most $(1+\varepsilon)t$. 
Then $w_\mu(t)\ls J(1+\varepsilon)t$. The upper bound follows from Theorem~\ref{th:near-minimal-measure-covering} and
the inequality $w_\mu(t)-w_\mu(a_nt)\ls w_\mu(t)$.
\end{proof}

\begin{corollary}\label{cor:depth-polytope-approximation}
Let $\mu$ be a log-concave probability measure on $\mathbb{R}^n$. For any $0<t<1/(4e)$ and $0<\varepsilon\ls1$, there exists a
polytope $P_{t,\varepsilon}$ such that
$$K_{(1+\varepsilon)t}\subseteq P_{t,\varepsilon}\subseteq\operatorname{int}(K_t),$$
with at most
$$B_n\varepsilon^{-(n-1)/2}\frac{w_\mu(t)-w_\mu(a_nt)}{t}$$
facets. Moreover,
$$\mu(K_t\setminus P_{t,\varepsilon})\ls w_\mu((1+\varepsilon)t)-w_\mu(t).$$
\end{corollary}

\begin{proof}We consider the covering half-spaces from Theorem~\ref{th:near-minimal-measure-covering} and set
$P_{t,\varepsilon}=\bigcap_j\big(\mathbb{R}^n\setminus\operatorname{int}(C_j)\big)$.
From the proof of Theorem~\ref{th:near-minimal-measure-covering} we see that $P_{t,\varepsilon}\subseteq\operatorname{int}(K_t)$.
Each covering half-space has measure at most $\sqrt{1+\varepsilon}\,t<(1+\varepsilon)t$,
so it cannot contain any point of $K_{(1+\varepsilon)t}$. This shows that $K_{(1+\varepsilon)t}\subseteq P_{t,\varepsilon}$. 
It also has non-empty interior, since $(1+\varepsilon)t<1/(2e)<q_*$. The set $P_{t,\varepsilon}$ is a bounded
intersection of at most $N$ closed half-spaces, hence a full-dimensional polytope with at most $N$ facets. Finally,
$K_t\setminus P_{t,\varepsilon}\subseteq K_t\setminus K_{(1+\varepsilon)t}$, and  taking measures we obtain the last 
inequality from the definition of $w_\mu$.
\end{proof}

In the last part of this section we verify that the exponent $(n-1)/2$ for the dependence on $\varepsilon$ in
Theorem~\ref{th:near-minimal-measure-covering} is optimal, and compute the associated shell complexity for the Gaussian measure.

Let $\gamma_n$ denote the standard Gaussian probability measure on $\mathbb{R}^n$.
Write $\phi(r)=(2\pi)^{-1/2}e^{-r^2/2}$,
$\overline\Phi(r)=\int_r^\infty\phi(s)\,ds$, and
$r_t=\overline\Phi^{-1}(t)$ for $0<t<1/2$. The closed geodesic cap with center at $u \in S^{n-1}$
and angular radius $\theta$ is the set $$C(u,\theta)=\{v \in S^{n-1}: \langle v,u\rangle \gr \cos \theta\}.$$

\begin{lemma}\label{lem:spherical-cover}Let $n\gr 2$ and write $\mathcal N_n(\theta)$ for the minimum number of closed
geodesic caps of angular radius $\theta$ whose union covers $S^{n-1}$. There exist $c_n,C_n>0$ such that
$$c_n\theta^{-(n-1)}\ls\mathcal N_n(\theta)\ls C_n\theta^{-(n-1)},\quad 0<\theta\ls1.$$
\end{lemma}

\begin{proof}Let $n\gr 3$. We integrate $(\sin s)^{n-2}$ on $[0,\theta]$, and using the fact that $\sin s\approx s$ in this interval, we see
that the surface measure of a cap of radius $\theta\ls1$ is comparable to $\theta^{n-1}$ up to constants depending only on the dimension. 
The case $n=2$ follows from the arc-length formula. Subadditivity gives the lower bound. 
Next, we choose a maximal set of centers with pairwise angular distances greater than $\theta$.
The caps of radius $\theta/2$ have disjoint interiors, and hence using the same area
estimate we see that the number of centers is bounded by $C_n\theta^{-(n-1)}$. By maximality,
the caps of radius $\theta$ cover the sphere, and we obtain the upper bound.
\end{proof}

\begin{theorem}\label{th:gaussian-covering}For every $n\gr2$ there exist constants $c_n,C_n,t_0>0$ such that
$$c_n\varepsilon^{-(n-1)/2}(\log(1/t))^{(n-1)/2}\ls N_{\gamma_n}(t,\varepsilon)\ls C_n\varepsilon^{-(n-1)/2}(\log(1/t))^{(n-1)/2}$$
for every $0<t<t_0$ and $0<\varepsilon\ls 1$. Moreover,
$$w_{\gamma_n}(t)/t\approx_n(\log(1/t))^{(n-1)/2},$$
and for every $a\in(0,1)$,
$$\frac{w_{\gamma_n}(t)-w_{\gamma_n}(at)}{t}\approx_{n,a}(\log(1/t))^{(n-1)/2}\quad(t\to 0^+).$$
\end{theorem}

\begin{proof}Let $X\sim \gamma_n$. For every $u\in S^{n-1}$ we have that $\langle u,X\rangle$ is a standard normal random
variable. Minimizing over all $u$ we get
$$q_{\gamma_n}(x)=\overline\Phi(\|x\|_2),\quad K_t=r_tB_2^n.$$
We set $b=r_t$, $a_t=r_{(1+\varepsilon)t}$ and $t_0=\overline\Phi(2)/2$. Note that $b\gr a_t>2$.
Every half-space of Gaussian measure at most $(1+\varepsilon)t$ can be written in the form 
$\{x:\langle u,x\rangle\gr h\}$ with $\|u\|_2=1$ and $h\gr a_t$.
We consider the larger half-space $\{x:\langle u,x\rangle\gr a_t\}$ whose measure is also at most
$(1+\varepsilon)t$. Its intersection with $bS^{n-1}$ is a spherical cap of angular radius
$$\theta=\arccos(a_t/b).$$
For any collection of such half-spaces we see that it forms a cover of $W_t=\{\|x\|_2 \gr b\}$ if and only if the intersections of these
half-spaces with $bS^{n-1}$ cover $bS^{n-1}$: this follows from the fact that along each outward ray, the chosen
positive inner product only increases. It follows that 
\begin{equation}\label{eq:gaussian-exact-spherical}
N_{\gamma_n}(t,\varepsilon)=\mathcal N_n\bigl(\arccos(a_t/b)\bigr).
\end{equation}
Integration by parts shows that
$$\frac{r}{r^2+1}\phi(r)\ls\overline\Phi(r)\ls\frac{\phi(r)}{r},\qquad r>0.$$
Therefore, the function $h(r)=\phi(r)/\overline\Phi(r)$ satisfies $r\ls h(r)\ls r+1/r$. Set $\eta=\log(1+\varepsilon)$. Then
$$\eta=\int_{a_t}^b h(r)\,dr,\quad\frac{b^2-a_t^2}{2}\ls\eta\ls\log2.$$
Since $a_t>2$, this implies that $b/a_t\ls\sqrt{1+(\log2)/2}$. Thus, on $[a_t,b]$ we have that $h(r)$ is comparable to $b$,
up to absolute constants. It follows that
$$b-a_t\approx\eta/b,\quad 1-a_t/b\approx\eta/b^2.$$
The same ratio bound shows that $0<\theta<1$. Indeed, from $1-\cos\theta\approx\theta^2$ and the fact that 
$\varepsilon/2\ls\eta\ls\varepsilon$ we get
$$\theta\approx\frac{\sqrt\varepsilon}{b}.$$
The tail bounds also give $b\sim\sqrt{2\log(1/t)}$ as $t\to 0^+$. Applying Lemma~\ref{lem:spherical-cover} to
\eqref{eq:gaussian-exact-spherical} we obtain the covering estimate.

Finally, polar integration gives
$$w_{\gamma_n}(t)=d_n\int_b^\infty r^{n-1}e^{-r^2/2}\,dr\sim d_nb^{n-2}e^{-b^2/2},$$
where $d_n>0$ depends only on $n$. Since $b\sim\sqrt{2\log(1/t)}$, we see that there exists a constant $d'_n>0$ such that
$$w_{\gamma_n}(t)\sim d'_nt(\log(1/t))^{(n-1)/2} \quad  (t\to 0^+) $$
For fixed $a\in(0,1)$, applying the last equivalence at $at$ and at $t$ we get
$w_{\gamma_n}(at)/w_{\gamma_n}(t)\to a$, which proves the shell estimate.
\end{proof}

\begin{remark}\rm
The Gaussian sharpness argument is independent of the general near-minimal covering theorem.
For any fixed sufficiently small $t$, it gives a positive multiple of
$\varepsilon^{-(n-1)/2}$ as a lower bound when $\varepsilon\to 0^+$.
Therefore no smaller exponent of $\varepsilon^{-1}$ can replace $(n-1)/2$ in
Theorem~\ref{th:near-minimal-measure-covering}, even if the shell factor is
replaced by $w_\mu(t)/t$. By affine invariance the same conclusions hold
for every nondegenerate Gaussian measure. In dimension one, two tails
suffice and are necessary for all sufficiently small $t$, so the corresponding
exponent is zero.
\end{remark}

\section{Optimal containment coverings for small half-spaces}\label{sec:containment}

A covering of the wet part by half-spaces of small measure does not, by itself, imply that each entire small half-space lies in one member of
the covering. Nevertheless, the whole half-space containment problem has an economical solution (Theorem~\ref{th:containment-family-optimal} and Corollary~\ref{cor:containment-optimal}).

Throughout this section, if $\mu$ is a log-concave probability measure on $\mathbb{R}^n$, we denote by $f$ its log-concave and upper-semicontinuous density. We set $d=n+1$ and consider the function $g=f^{-1/d}$, extending it by $+\infty$ outside the set $\{f>0\}$. 

\begin{lemma}\label{lem:containment-projective-lift}
Let $\mu$ be a log-concave probability measure on $\mathbb{R}^n$. Then $g$ is convex and
\begin{equation}\label{eq:containment-perspective-body}
\mathcal K_f=\overline{\left\{(y,r)\in\mathbb{R}^n\times(0,\infty):r g(y/r)\ls1\right\}}
\end{equation}
is a  convex body in $\mathbb{R}^d$ with $0 \in \mathcal K_f$ and $\operatorname{vol}_d(\mathcal K_f)=1/d$. Moreover, 
for every Borel set $B \subseteq \mathbb{R}^n$ we have \begin{equation}\label{eq:containment-projective-pushforward}
\mu(B)=\frac{\operatorname{vol}_d\bigl(\{(y,r)\in\mathcal K_f:r>0,\ y/r\in B\}\bigr)}{\operatorname{vol}_d(\mathcal{K}_f)}.
\end{equation} 
\end{lemma}

\begin{proof}
Since $f$ is log-concave and upper-semicontinuous, the function $g=\exp((-\log f)/d)$ is convex and lower-semicontinuous. 
Thus the map $(y,r)\mapsto rg(y/r)$ is convex and lower-semicontinuous. Its sublevel set in \eqref{eq:containment-perspective-body} 
is convex and relatively closed in $\mathbb{R}^n \times (0,\infty)$. It follows that its closure $\mathcal K_f$ is convex and 
adds points only on $\mathbb{R}^n \times \{0\}$, a set of $d$-dimensional Lebesgue measure zero.

Choose $x_0$ with $f(x_0)>0$. Then $(rx_0,r)$ belongs to the defining set for every $0<r\ls f(x_0)^{1/d}$, and letting $r\to0^+$ 
we see that $0\in\mathcal K_f$. 

For a Borel set $B\subseteq\mathbb R^n$, Tonelli's theorem gives
$$\mu(B)=d\int_B\int_0^{f(x)^{1/d}}r^n\,dr\,dx.$$
The Jacobian of $(x,r)\mapsto(y,r)=(rx,r)$ is $r^n$. Since taking the closure we do not add points with $r>0$, we obtain
$$\frac{1}{d}\mu(B)=\operatorname{vol}_d\bigl(\{(y,r)\in\mathcal K_f:r>0,\ y/r\in B\}\bigr).$$
Choosing $B=\mathbb R^n$ we get $\operatorname{vol}_d(\mathcal K_f)=1/d$. Thus, the preceding identity also proves \eqref{eq:containment-projective-pushforward}. Finally, since $\mathcal{K}_f$ is closed, convex and has finite and positive volume, 
it is a convex body. 
\end{proof}

Every closed half-space
$H=\{x:\ip{u}{x}\gr a\}$ corresponds to the cap $$C_H=\mathcal K_f\cap\{(y,r):\ip{u}{y}\gr ar\}.$$ Lemma~\ref{lem:containment-projective-lift} then yields
\begin{equation}\label{eq:containment-projective-cap}
\mu(H)= \frac{\operatorname{vol}_d(C_H)}{\operatorname{vol}_d(\mathcal K_f)}.
\end{equation}

\begin{lemma}\label{lem:radial-transfer}
Let $R\subseteq\mathbb{R}^{n+1}$ be a closed convex set containing the origin and set
$$r_x=f(x)^{1/(n+1)},\quad p(x)=(r_xx,r_x)\quad(f(x)>0),$$
$$A(R)=\{x:f(x)>0,\ p(x)\in R\}\cup\{x:f(x)=0\}.$$
Then $A(R)$ is a Borel set and
\begin{equation}\label{eq:radial-transfer-bound}
\mu(A(R))\ls\frac{\operatorname{vol}_{n+1}(\mathcal K_f\cap R)}
{\operatorname{vol}_{n+1}(\mathcal K_f)}.
\end{equation}
In addition, if $C_H\subseteq R$ for the cap associated to a closed half-space $H$,
then $H\subseteq A(R)$.
\end{lemma}

\begin{proof}
The density $f$ is Borel measurable, so $A(R)$ is a Borel set. The set $\{f=0\}$ has zero $\mu$-measure. By Lemma~\ref{lem:containment-projective-lift}, $\mu(A(R))$ is equal to the normalized volume of the points $(y,r)\in\mathcal K_f$ with $r>0$, $f(y/r)>0$, and $p(y/r)\in R$.

For any such point, set $x=y/r$. Then $0<rg(x)\ls1$ and
$$
p(x)=\left(\frac{x}{g(x)},\frac{1}{g(x)}\right),\qquad
(y,r)=rg(x)\,p(x).
$$
Since $R$ is convex and contains both $0$ and $p(x)$, we have $(y,r)\in R$. Comparing volumes proves \eqref{eq:radial-transfer-bound}.

Finally, if $C_H\subseteq R$ and $x\in H$ with $f(x)>0$, then $p(x)\in C_H\subseteq R$, so $x\in A(R)$. Points of $H$ with $f(x)=0$ belong to $A(R)$ by definition. Thus $H\subseteq A(R)$.
\end{proof}

We will also use the following result. It is a cap-covering lemma of Arya, da Fonseca, and
Mount~\cite[Lemma~4.3]{Arya-daFonseca-Mount-2026}.

\begin{lemma}\label{lem:restricted-family-classical-cap-covering}
For $d\gr2$ let $K\subseteq\mathbb{R}^d$ be a convex body, and $\mathcal C$ be any family of positive-volume caps
$C=K\cap H$. There exist pairwise disjoint convex bodies
$M_i\subseteq K$, corresponding convex bodies $R_i\subseteq\mathbb{R}^d$ with
$M_i\subseteq R_i$ and a constant $c_d>0$, depending only on $d$, such that
$$\operatorname{vol}_d(R_i)\ls c_d\operatorname{vol}_d(M_i)$$
and for every $C\in\mathcal C$ there is an index $i$ with
$$M_i\subseteq C\subseteq R_i.$$
\end{lemma}

\begin{theorem}\label{th:containment-family-optimal}
For every $n\gr1$ there exists a constant $C_n>0$ with the following property. Let $\mu$ be a log-concave probability
measure on $\mathbb{R}^n$, $0<t<q_\ast$, and $\mathcal H$ any family of closed half-spaces, all of $\mu$-measure $t$. If $\,\bigcup_{H\in\mathcal H}H$ is a Borel set, then there exist
Borel sets $A_1,\ldots,A_J$ with
\begin{equation}\label{eq:containment-family-bounds}
J\ls C_n\frac{\mu\bigl(\bigcup_{H\in\mathcal H}H\bigr)}{t}, \quad\mu(A_i)\ls C_nt,
\end{equation}
such that every $H\in\mathcal H$ is contained in some $A_i$.
\end{theorem}

\begin{proof}
Set $K=\mathcal K_f$ and lift each $H\in\mathcal H$ to the cap $C_H$. Apply Lemma~\ref{lem:restricted-family-classical-cap-covering} in dimension $d=n+1$ to obtain pairs $(M_i,R_i)$. Retain precisely the indices
$$I:=\{i:\text{there exists }H\in\mathcal H\text{ with }M_i\subseteq C_H\subseteq R_i\}.$$
For each $i \in I$ choose a half-space $H_i\in\mathcal H$ with $M_i\subseteq C_{H_i}\subseteq R_i$. 
Using \eqref{eq:containment-projective-cap} and Lemma~\ref{lem:restricted-family-classical-cap-covering} we get
$$\frac{t}{c_d}\operatorname{vol}_d(K)\ls\operatorname{vol}_d(M_i)\ls t\operatorname{vol}_d(K),
\qquad \operatorname{vol}_d(R_i)\ls c_dt\operatorname{vol}_d(K).$$
Since they are pairwise disjoint subsets of $K$, there are finitely many sets $M_i$. Relabel them as $M_1,\ldots,M_J$. Using \eqref{eq:containment-projective-pushforward}, we see that
$$J\frac{t}{c_d}\operatorname{vol}_d(K)\ls\sum_{i=1}^J\operatorname{vol}_d(M_i)
\ls\operatorname{vol}_d\Big(\bigcup_{H\in\mathcal H}C_H\Big)=\operatorname{vol}_d(K)\mu\Big(\bigcup_{H\in\mathcal H}H\Big).$$
This proves the cardinality estimate.

Every cap $C_H$ contains the origin, because $0 \in K$. In particular, $C_{H_i}\subseteq R_i$ implies $0\in R_i$ for every $i \in I$. Define $A_i=A(R_i)$ using Lemma~\ref{lem:radial-transfer}. Then apply \eqref{eq:radial-transfer-bound}:
$$
\mu(A_i)\ls\frac{\operatorname{vol}_d(K\cap R_i)}{\operatorname{vol}_d(K)}\ls c_dt.
$$
Finally, let $H\in\mathcal H$. There exists an index $1 \ls i \ls J$ such that $M_i\subseteq C_H\subseteq R_i$. Lemma~\ref{lem:radial-transfer} then gives $H\subseteq A(R_i)=A_i$. This proves the containment for every $H\in\mathcal H$.
\end{proof}

\begin{remark}\rm
The cardinality bound is optimal up to dimensional constants. Indeed, any finite family of Borel sets $A_1,\ldots,A_J$ with $\mu(A_i)\ls C_nt$ that contains every member of $\mathcal H$ in some $A_i$ satisfies
$$
\mu\Bigl(\bigcup_{H\in\mathcal H}H\Bigr)\ls\sum_{i=1}^J\mu(A_i)\ls JC_nt.
$$
\end{remark}

\begin{remark}\rm In the proof of Theorem~\ref{th:containment-family-optimal} we only use the fact that $g=f^{-1/(n+1)}$ is convex 
and that the perspective body in \eqref{eq:containment-perspective-body} is compact. For such densities, Theorem~\ref{th:containment-family-optimal} 
continues to hold, with the same constants. An example is given by the standard Cauchy density $f(x)\propto(1+\|x\|_2^2)^{-(n+1)/2}$:
we have that $g(x)$ is proportional to $\sqrt{1+\|x\|_2^2}$ and its perspective body is the upper half of a Euclidean ball.
\end{remark}

\begin{corollary}\label{cor:containment-optimal}Let $C_n$ be the constant from Theorem~$\ref{th:containment-family-optimal}$.
For every log-concave probability measure $\mu$ on $\mathbb{R}^n$ and any $0<t<q_\ast$,
there exist Borel sets $A_1,\ldots,A_J$ such that
\begin{equation}\label{eq:containment-optimal}
J\ls C_n\frac{w_\mu(t)}{t},\quad\mu(A_i)\ls C_nt,
\end{equation}
and every closed half-space with $\mu$-measure at most $t$ is contained in at least one $A_i$.
\end{corollary}

\begin{proof}We apply Theorem~\ref{th:containment-family-optimal} to the family of all closed half-spaces of $\mu$-measure exactly $t$, because every half-space of smaller measure can be enlarged to a half-space of measure $t$, by continuity of one-dimensional marginals. The union of these half-spaces equals $\{q_\mu\ls t\}$, which is a Borel set and has measure $w_\mu(t)$.
\end{proof}

\begin{corollary}\label{cor:random-hyperplane-bound}Let $n\gr2$ and let $Y_1,\ldots,Y_n$ be i.i.d random vectors with a log-concave probability distribution $\mu$ on $\mathbb{R}^n$. Almost surely their affine hull is a hyperplane $H_Y$. Write
$H_Y^+,H_Y^-$ for its two closed sides and define
$$G_\mu(t)=\mathbb{P}\bigl(\min\{\mu(H_Y^+),\mu(H_Y^-)\}\ls t\bigr).$$
There exists a constant $c_n>0$, depending only on $n$, such that
\begin{equation}\label{eq:random-hyperplane-bound}
G_\mu(t)\ls c_n t^{n-1}w_\mu(t),\qquad 0<t<q_\ast.
\end{equation}
\end{corollary}

\begin{proof}Fix $t$ in the range of Corollary~\ref{cor:containment-optimal}, and let $A_1,\ldots,A_J$ be the Borel sets in the
statement of the corollary. If
$\min\{\mu(H_Y^+),\mu(H_Y^-)\}\ls t$, one of the two closed sides of $H_Y$ is contained in some $A_i$. Since all $Y_j$ lie on $H_Y$,
this implies $Y_1,\ldots,Y_n\in A_i$. It follows that
$$G_\mu(t)\ls\sum_{i=1}^J\mathbb{P}(Y_1,\ldots,Y_n\in A_i)=\sum_{i=1}^J\mu(A_i)^n\ls J(C_nt)^n\ls C_n^{n+1} t^{n-1} w_\mu(t),$$
where $C_n$ is the constant in \eqref{eq:containment-optimal}.
\end{proof}

\section{Wet-part estimates}\label{sec:shell}

The classical wet-part estimate for convex bodies of B\'ar\'any--Larman~\cite{Barany-Larman-1988}, as presented by B\'ar\'any in~\cite[Corollary~11.4]{Barany-2008-notes}, has sharp exponent $2/(n+1)$. We first state the exact one-dimensional formula and transfer this absolute bound to log-concave measures by layer-cake integration. We then combine that transfer with the economic covering theorem to prove a relative comparison between depth levels, with consequences for depth moments and covering numbers.

\begin{proposition}\label{prop:one-dimensional-wet}Let $\mu$ be a log-concave probability measure on $\mathbb{R}$. Then, for $0<t<q_\ast$,
\begin{equation}\label{eq:one-dimensional-wet}
w_\mu(t)=2t.
\end{equation}
\end{proposition}

\begin{proof}Let $F$ be the distribution function of $\mu$.  Since $\mu$ has a density, $F$ is continuous, and in dimension one
$$q_\mu(x)=\min\{F(x),1-F(x)\}.$$
Thus $W_t$ is the union of the two tails of measure $t$, which proves \eqref{eq:one-dimensional-wet}.
\end{proof}

For a convex body $K$ in $\mathbb{R}^n$, to distinguish cap volume from the normalized depth $v_K$ used in the introduction, write
$$\operatorname{Wet}_K(\varepsilon)
:=\bigl\{x\in K:\inf_{H\ni x} \operatorname{vol}_n(K\cap H)\ls \varepsilon\bigr\}
=U_{\varepsilon/\operatorname{vol}_n(K)}(K).$$

\begin{lemma}\label{lem:classical-wet-upper}For $n\gr2$, let $\overline{b}_n$ be the least constant such that, for every
convex body $K$ in $\mathbb{R}^n$ and every $\varepsilon>0$,
\begin{equation}\label{eq:classical-wet-upper}
\operatorname{vol}_n\bigl(\operatorname{Wet}_K(\varepsilon)\bigr) \ls \overline{b}_n \operatorname{vol}_n(K)^{\frac{n-1}{n+1}}\varepsilon^{\frac{2}{n+1}}.
\end{equation}
Then $\overline{b}_n<\infty$.
\end{lemma}

\begin{proof}
The classical convex-body result that was mentioned above (see~\cite{Barany-Larman-1988,Barany-2008-notes}) gives
constants $c_n>0$ and $\varepsilon_n>0$, depending only on $n$, such that every convex body $K$ of volume $1$
in $\mathbb{R}^n$ satisfies
$$\operatorname{vol}_n\bigl(\operatorname{Wet}_K(\varepsilon)\bigr)\ls c_n\varepsilon^{2/(n+1)}
\quad(0<\varepsilon<\varepsilon_n).$$
For $\varepsilon\gr\varepsilon_n$, the bound by $1$ extends this estimate
with constant $\widetilde c_n=\max\{c_n,\varepsilon_n^{-2/(n+1)}\}$.
For general volume $V=\operatorname{vol}_n(K)$, we apply the previous estimate to $V^{-1/n}K$:
$$\operatorname{vol}_n\bigl(\operatorname{Wet}_K(\varepsilon)\bigr)
=V\operatorname{vol}_n\bigl(\operatorname{Wet}_{V^{-1/n}K}(\varepsilon/V)\bigr)
\ls \widetilde c_nV^{(n-1)/(n+1)}\varepsilon^{2/(n+1)}.$$
Hence, $\overline{b}_n\ls \widetilde c_n<\infty$.
\end{proof}

\begin{theorem}\label{th:universal-wet-upper}For every $n\gr1$ there exists $D_n>0$ such that every log-concave probability measure $\mu$ on $\mathbb{R}^n$ satisfies
\begin{equation}\label{eq:universal-wet-upper}
 w_\mu(t)\ls D_n t^{\frac{2}{n+1}}, \quad t>0.
\end{equation}
The exponent $2/(n+1)$ cannot be increased uniformly over all such measures.
\end{theorem}

\begin{proof}The case $n=1$ follows from Proposition~\ref{prop:one-dimensional-wet}. Assume that $n\gr2$ and set
$$\alpha:=\frac{2}{n+1}.$$
The statement is affinely invariant, so we may assume that the density $f$ of $\mu$ satisfies $\|f\|_\infty=1$.  For $0<s<1$, set
$$L_s:=\{x:f(x)\gr s\}.$$
These sets are bounded convex bodies, and
\begin{equation}\label{eq:layer-cake-normalization}
 1=\int_0^1 \operatorname{vol}_n(L_s)\,ds.
\end{equation}
Fix $x\in W_t\cap L_s$.  There exists a half-space $H$ containing $x$ with $\mu(H)\ls t$.  
Since $L_s\subseteq L_r$ for $0<r<s$, we have
$$s\operatorname{vol}_n(H\cap L_s)\ls\int_0^s \operatorname{vol}_n(H\cap L_r)\,dr\ls\mu(H)\ls t.$$
It follows that
\begin{equation}\label{eq:level-wet-inclusion}
 W_t\cap L_s\subseteq \operatorname{Wet}_{L_s}(t/s).
\end{equation}
Using Lemma~\ref{lem:classical-wet-upper}, we may write
\begin{equation}\label{eq:wet-level-integral}
 w_\mu(t) \ls \overline{b}_n t^\alpha \int_0^1 s^{-\alpha}\operatorname{vol}_n(L_s)^{1-\alpha}\,ds,
\end{equation}
and it remains to bound the last integral.  Set
$$\beta:=1-\alpha=\frac{n-1}{n+1},\quad\varphi(r):=\operatorname{vol}_n(L_{e^{-r}})^{1/n}.$$
The function $\varphi$ is non-decreasing and concave by the
log-concavity of $f$ and the Brunn-Minkowski inequality.  Moreover,
\begin{equation}\label{eq:level-volume-normalization}
\int_0^\infty e^{-r}\varphi(r)^n\,dr=1.
\end{equation}
Since $\varphi$ is concave, we have
$$ \varphi(u/\beta)\ls\beta^{-1}\varphi(u),\quad u\gr0.$$
It follows that
\begin{align*}
\int_0^1s^{-\alpha}|L_s|^{1-\alpha}\,ds &=\int_0^\infty e^{-\beta r}\varphi(r)^{n\beta}\,dr\\
&\ls \beta^{-1-n\beta} \int_0^\infty e^{-u}\varphi(u)^{n\beta}\,du\\
&\ls \beta^{-1-n\beta}.
\end{align*}
For the last inequality we use \eqref{eq:level-volume-normalization} and the monotonicity of $L_p$ norms for the probability 
measure $e^{-u}\,du$. Choosing
\begin{equation}\label{eq:wet-constant-transfer}
 D_n=\overline{b}_n\,\beta^{-1-n\beta},
\end{equation}
we obtain \eqref{eq:universal-wet-upper}. In the case of the Euclidean ball we have $w(t)\approx_n t^{2/(n+1)}$, and this
shows that the exponent is optimal.
\end{proof}

\begin{corollary}\label{cor:wet-constant-comparison}Let $n\gr2$ and write $\overline{D}_n$ for the smallest constant such that
$$w_\mu(t)\ls \overline{D}_nt^{2/(n+1)}$$
for every log-concave probability measure $\mu$ on $\mathbb{R}^n$ and every $t>0$. Then, 
\begin{equation}\label{eq:wet-constant-comparison}
\overline{b}_n\ls \overline{D}_n \ls \beta^{-1-n\beta}\overline{b}_n,
\end{equation}
where $\beta=(n-1)/(n+1)$. In particular, $\overline{b}_n$ and $\overline{D}_n$ are comparable up to absolute constants.
\end{corollary}

\begin{proof}The upper bound follows from \eqref{eq:wet-constant-transfer}. For the lower bound, we apply Theorem~\ref{th:universal-wet-upper} 
to the uniform measure on a convex body. Note that 
$$(1+n\beta)\log(1/\beta)\ls \frac{2n}{n-1}\ls 4, $$
and hence $\beta^{-1-n\beta}$ is bounded uniformly in $n\gr2$.
\end{proof}

\begin{theorem}\label{th:optimal-wet-relative}For every integer $n\gr1$ there exist constants $B_n>0$ and $t_n>0$ such that, for
every log-concave probability measure $\mu$ on $\mathbb{R}^n$,
\begin{equation}\label{eq:optimal-wet-relative}
 w_\mu(s) \ls B_n\left(\frac{s}{t}\right)^{\frac{2}{n+1}}w_\mu(t), \quad 0<s\ls t<t_n.
\end{equation}
The exponent is optimal.
\end{theorem}

\begin{proof}We set $\alpha=2/(n+1)$ and apply Theorem~\ref{th:economic-cap-covering} to obtain half-spaces $C_1,\ldots,C_N$ such that
$$W_t\subseteq\bigcup_{i=1}^N C_i,\quad m_i:=\mu(C_i)\ls A_nt,\quad N\ls D_n\frac{w_\mu(t)}{t}.$$
If $m_i>0$ then we define
$$\nu_i=\frac{\mu(\,\cdot\cap C_i)}{m_i}.$$
Note that $\nu_i$ is a log-concave probability measure.  If $x\in W_s\cap C_i$, then a half-space $H$ containing $x$ with $\mu(H) \ls s$ satisfies
$$\nu_i(H)\ls\frac{\mu(H)}{m_i}\ls\frac{s}{m_i}.$$
It follows that
$$W_s\cap C_i\subseteq W_{s/m_i}^{\nu_i}.$$
Then Theorem~\ref{th:universal-wet-upper} shows that
$$\mu(W_s\cap C_i)\ls D_n^\prime s^\alpha m_i^{1-\alpha}\ls D_n^\prime A_n^{1-\alpha}s^\alpha t^{1-\alpha}.$$
The last inequality also holds for those $i$ with $m_i=0$. Since $W_s\subseteq W_t$, summing over the covering half-spaces we get
$$w_\mu(s) \ls D_n^\prime A_n^{1-\alpha}N s^\alpha t^{1-\alpha}\ls B_n \left(\frac{s}{t}\right)^\alpha w_\mu(t).$$
The uniform measure on the Euclidean ball again shows that the exponent in Theorem~\ref{th:optimal-wet-relative} cannot be increased.
\end{proof}

\begin{corollary}\label{cor:depth-negative-moments}Let $\mu$ be a log-concave probability measure on $\mathbb{R}^n$. If $\alpha=2/(n+1)$ and $0<p<\alpha$, then
\begin{equation}\label{eq:depth-negative-moment}
 \mathbb{E}_\mu \big[q_\mu^{-p}\big]\ls 1+\frac{pD_n}{\alpha-p},
\end{equation}
for some dimensional constant $D_n>0$. The range $p<2/(n+1)$ is optimal uniformly over log-concave probability measures.
\end{corollary}

\begin{proof} The absolute wet-part bound in Theorem~\ref{th:universal-wet-upper} gives $w_\mu(r)\ls D_nr^\alpha$.  Since $q_\mu\ls1$, integration of the distribution function gives
$$\mathbb{E}_\mu \big[q_\mu^{-p}\big]=1+p\int_1^\infty u^{p-1} w_\mu(u^{-1})\,du \ls1+\frac{pD_n}{\alpha-p}.$$
For the uniform measure on the Euclidean ball, $w_\mu(r)\approx_n r^\alpha$, so $\mathbb{E}_\mu \big[q_\mu^{-p}\big]=\infty$ for $p\gr\alpha$.
\end{proof}

\begin{corollary}\label{cor:universal-covering-cardinality}For every log-concave probability measure $\mu$ on $\mathbb{R}^n$ and every $0<t<t_n$, the covering half-spaces in Theorem~$\ref{th:economic-cap-covering}$ may be chosen with
\begin{equation}\label{eq:universal-covering-cardinality}
N\ls D_n t^{-\frac{n-1}{n+1}}.
\end{equation}
The exponent is optimal over all log-concave probability measures.
\end{corollary}

\begin{proof}
The upper bound follows by combining Theorem~\ref{th:economic-cap-covering}
with Theorem~\ref{th:universal-wet-upper}. For sharpness, let $\mu$ be the
uniform probability measure on a Euclidean ball. Then
$w_\mu(t)\approx_n t^{2/(n+1)}$ as $t\to 0^+$. Any cover by $N$
half-spaces of measure at most $A_nt$ must satisfy
$$w_\mu(t)\ls\sum_{i=1}^N\mu(C_i)\ls NA_nt,$$
and hence
$$N\gr\frac{w_\mu(t)}{A_nt}\gr c_n t^{-(n-1)/(n+1)}$$
for all sufficiently small $t$.
The same argument in dimension one uses $w_\mu(t)=2t$ and gives a positive
constant lower bound, as required for exponent zero.
\end{proof}

\bigskip

\noindent {\bf Acknowledgements.} I would like to thank A.~Giannopoulos for helpful discussions.

\bigskip

\noindent\textbf{AI Usage Statement.} The author acknowledges the assistance of ChatGPT (OpenAI) at the final stage of the preparation of this work: this included proof checking, corrections and suggestions on the presentation.  The author takes full responsibility for the content and correctness of the manuscript.

\bigskip

\bigskip 

\medskip

\noindent\textbf{2020 Mathematics Subject Classification.} Primary 52A20, 52C17; Secondary 52A22, 52A27, 60D05.

\smallskip 

\noindent\textbf{Keywords.} Log-concave measures, economic cap coverings, Tukey half-space depth, wet parts, floating bodies, Macbeath regions, containment coverings, random polytopes.

\bigskip

\bigskip

\bigskip

\bigskip

\noindent \textsc{Minas \ Pafis}: Department of Mathematics, National and Kapodistrian University of Athens, Panepistimioupolis 157-84,
Athens, Greece.

\smallskip

\noindent \textit{E-mail:} \texttt{mipafis@math.uoa.gr}

\end{document}